\documentclass[a4paper, oneside, 11pt]{article}

\usepackage{a4wide}
\usepackage{graphicx}
\usepackage{pb-diagram, lamsarrow,pb-lams}
\usepackage{amssymb}
\usepackage{parskip}
\usepackage{theorem}
\usepackage{enumerate}
\usepackage{floatflt}
\usepackage[centertags]{amsmath}
\usepackage[latin1]{inputenc}

\input xy
\xyoption{all} \CompileMatrices

\renewcommand{\b}{\bullet}

\newcommand{\Hom}{\mbox{\upshape{Hom}}}

\newcommand{\ind}{\mbox{\upshape{ind}}}
\renewcommand{\ker}{\mbox{\upshape{ker}}}
\renewcommand{\dim}{\mbox{\upshape{dim}}}
\newcommand{\Det}{\mbox{\upshape{Det}}}
\newcommand{\Qu}{\mbox{\upshape{Quot}}}
\newcommand{\res}{\mbox{\upshape{res}}}
\newcommand{\cd}{\mbox{\upshape{cd}}}

\newcommand{\Q}{\mathbb{Q}}
\newcommand{\N}{\mathbb{N}}

\newcommand{\Z}{\mathbb{Z}}
\newcommand{\F}{\mathbb{F}}
\newcommand{\TT}{\overline{T}}
\newcommand{\Gal}{\mbox{\upshape{Gal}}}
\newcommand{\QQ}{\mathcal{Q}}

\renewcommand{\o}{\mathfrak{o}}

\newcommand{\m}{\mathfrak{m}}
\newcommand{\n}{\mathfrak{n}}

\newcommand{\Tr}{\mbox{\upshape{Tr}}}
\renewcommand{\mod}{\:\mbox{\upshape{mod}}\:}

\newcommand{\p}{\mathfrak{p}}
\newcommand{\ch}{\mbox{\upshape{char}}}

\newcommand{\NN}{\mbox{\upshape{N}}}

\newcommand{\nr}{\mbox{\upshape{nr}}}

\newcommand{\gal}{\mbox{\scriptsize{\upshape{Gal}}}}
\newcommand{\ur}{\mbox{\scriptsize{\upshape{nr}}}}

\newtheorem{lemma}{Lemma}
\newtheorem{theorem}{Theorem}
\newtheorem{proposition}{Proposition}
\newtheorem{corollary}{Corollary}
\newtheorem{remark}{Remark}
\newtheorem{definition}{Definition}

\begin{document}

\title{On the Iwasawa algebra for pro-$l$ Galois groups}

\author{Irene Lau}\date{}
\maketitle

\begin{abstract}
Let $l$ be an odd prime and $K/k$ a Galois extension of totally real number fields with Galois group $G$ such that $K/k_\infty$ and $k/\Q$ are finite. We give a full description of the algebraic
structure of the semisimple algebra $\QQ G=\Qu(\Lambda G)$ for pro-$l$ Galois groups $G$ with $\Lambda G = \Z_l[[G]]$ the Iwasawa algebra. We moreover compute the cohomological dimension of the centres of the Wedderburn components of $\QQ G$ and state some results on the  completion of $\QQ G$.
\end{abstract}

\section{Introduction}
\label{intro}
We fix an an odd prime number $l$ and a Galois extension $K/k$ of totally real number fields with Galois group $G$ such that $K/k_\infty$ and $k/\Q$ are finite. As usual, $k_\infty$  is the cyclotomic $\Z_l$-extension of $k$.
Next, the Iwasawa algebra $\Lambda G=\Z_l[[G]]$ denotes the completed group ring of $G$ over $\Z_l$ and $\QQ G=\Qu(\Z_l[[G]])$ is its total ring of fractions with respect to all central non-zero divisors. $\QQ G$ finds its way into non-commutative Iwasawa theory via the localization sequence of $K$-theory
$$\rightarrow K_1(\Lambda G)\rightarrow K_1(\QQ G)\stackrel{\partial}{\rightarrow}K_0T(\Lambda G)\rightarrow$$
and a determinant map
$$\Det: K_1(\QQ G)\rightarrow \Hom(R_lG, (\Q_l^c\otimes_{\Q_l}\QQ \Gamma_k)^\times).$$
As in the classical case of Iwasawa,  Ritter and Weiss link a $K$-theoretic substitute $\mho\in K_0T(\Lambda G)$ of the Iwasawa module $X$ to the Iwasawa $L$-function which is derived from the $S$-truncated Artin $L$-function for a finite set $S$ of places of $k$ containing all archimedian ones and those which ramify in $K$. This Iwasawa $L$-function lies in the upper $\Hom$-group. Here, as usual, $X=G(M/K)$ with $M$ denoting the maximal abelian $l$-extension of $K$ which is unramified outside $S$.

With this, the main conjecture of equivariant Iwasawa theory says
$$\exists!\ \Theta\in K_1(\QQ G):\ \Det(\Theta)=L\ \&\ \partial(\Theta)=\mho.$$
The uniqueness of $\Theta$ would follow from a conjecture by Suslin being true. This conjecture applied to our situation means  that $SK_1(\QQ G)=\ker(\Det)=1$ (we suppress here the condition on the cohomological dimension of the centres which we show  to be fulfilled in the last section).

Recently (compare \cite{Iw10}), Ritter and Weiss gave a complete proof of this main conjecture up to its uniqueness statement whenever Iwasawa's $\mu$-invariant vanishes. In \cite{Kak2}, Kakde also gave a proof. In fact, he does not restrict to 1-dimensional $l$-adic Lie groups as Ritter and Weiss do but gives a proof for higher dimensional admissible $l$-adic Lie extensions satisfying the conjecture $\mu=0$. Yet, the uniqueness stated and proven in his main conjecture is weaker than in the formulation of Ritter and Weiss because he does not consider the Whitehead group $K_1(\QQ G)$ but its quotient by the reduced Whitehead group $SK_1(\QQ G)$. Thus, the question whether $\Theta$ is unique in $K_1(\QQ G)$ is still open.

Aiming to prove the uniqueness of $\Theta$, the algebra $\QQ G$ is of decisive interest. Beyond that, it is an interesting algebraic object in itself and  the results achieved here may be interesting also from a purely algebraic point of view because they show that properties of group algebras of finite groups do generally not persist when passing to projective limits of such. 

In Section 2, we resolve the structure of the algebra $\QQ G$ for pro-$l$ Galois groups $G$. Here, we restrict ourselves to the case of pro-$l$ groups as this is the crucial case for the proof of the main conjecture. Roquette  (see \cite{Ro}) showed that for finite $l$-groups $H$, the group algebra $F[H]$ over a field $F$ of characteristic 0 splits into full matrix rings over fields, i.e.~that skew fields do not appear in the Wedderburn components of $F[H]$. Although $\QQ G=\Qu(\Z_l[[G]])$ is not a group algebra, one might expect
that there do not occur Schur indices in its structure because $G$ is as pro-$l$ group the projective limit of finite $l$-groups. However, this is not the case: As we show in Section 2, non-trivial Schur indices appear, but the occuring skew fields are all cyclic. We resolve the algebraic structure of the Wedderburn components of $\QQ G$ completely. Moreover, we give an example in which such a non-trivial case appears as Galois group in the equivariant Iwasawa setting.

In the third section, we consider fields of type $\Q_l(\zeta)\otimes_{\Q_l}\QQ \Gamma$ with $\zeta$ a primitive $l$-power root of unity and $\Gamma\cong \Z_l$. We compute the cohomological dimensions of these fields and their $\p$-adic completions for primes $\p$ of height 1. Furthermore, we use this to show that the Suslin conjecture  is true for the completed $\QQ_\wedge G$ if the completion is not with respect to the prime above $l$. For the completion with respect to $l$, we are at least able to compute the completion of the skew fields underlying the simple components of $\QQ G$ and their residue skew fields. Furthermore, we compute the cohomological dimension of the centres of the Wedderburn components of $\QQ G$.
We  hope that these results might be useful for future work on Suslin's conjecture.

In this paper, I recollect some results of my PhD thesis.
I would like to thank my supervisor Jürgen Ritter who encouraged me to work on this fascinating subject. He was and still is a valuable mentor.

\section{The structure of $\QQ G$}
\label{sec:1}
In this section, we restrict ourselves to pro-$l$ Galois groups $G$ and compute the structure of $\QQ G$ for these groups.

We again fix an an odd prime number $l$ and a Galois extension $K/k$ of totally real number fields with Galois group $G$ such that $k/\Q$ and $K/k_\infty$ are finite. As usual, $k_\infty$  denotes the cyclotomic $\Z_l$-extension of $k$ and thus $G(k_\infty/k)\cong \Z_l$. Moreover, we set $H:=G(K/k_\infty)$ which is a finite $l$-group by the above. Recall that $G=H\rtimes \Gamma$ splits with $\Gamma \cong \Z_l$ (see e.g. \cite[p. 551]{Iw2}).

 Next, the Iwasawa algebra $\Lambda G=\Z_l[[G]]$ denotes the completed group ring of $G$ over $\Z_l$ and $\QQ G=\Qu(\Z_l[[G]])$ is its total ring of fractions with respect to all central non-zero divisors. There exists an $m\in \N$ such that $\Gamma_0=\Gamma^{l^m}$ is a central subgroup of $G$ and with this
$$\QQ G=\bigoplus_{i=0}^{l^m-1}(\QQ \Gamma_0)[H]\gamma^i.$$

 Let $A$ be a simple component of the finite dimensional semisimple $\QQ \Gamma$-algebra $\QQ G$. Then, $A$ corresponds to an irreducible $\chi\in R_l G$ (modulo $W$-twists and the action of the Galois group $G(\Q_l^c/\Q_l)$) by $A e_\chi \neq 0$. 
From now on, we fix a representative $\chi$ in the orbit under $W$-twisting.

For the investigation of the structure of $A$, we need  the following 
property:

\begin{lemma}
Let $A$ be a simple component of $\QQ G=\bigoplus_{i=0}^{l^m-1}(\QQ\Gamma_0)[H]\gamma^{i}$. Then
$$A\cap(\QQ\Gamma_0)[H]\neq 0.$$
\end{lemma}
\textbf{Proof:}
First, we know that there exists a primitive central idempotent $e\in \QQ G$ with $A=\QQ Ge$. Furthermore, the tensor product $\Q_l^c\otimes_{\Q_l}A=\Q_l^c\otimes_{\Q_l}\QQ Ge$ is the direct sum of some simple components which correspond to central primitive  idempotents $e_{\chi_\nu}\in \QQ^cG$  for a finite set of irreducible $\chi_\nu\in R_lG$. Thus, we get
$$\Q_l^c\otimes_{\Q_l}A=\Q_l^c\otimes_{\Q_l}\QQ Ge=\bigoplus_\nu{\QQ^cGe_{\chi_\nu}}.$$
As the $e_{\chi_\nu}$ are exactly the primitive central idempotents lying over $e$, we conclude $e=\sum_{\nu} e_{\chi_\nu}\in A$. Therefore, by $$e_{\chi_\nu}=\sum_{\eta\mid \mbox{\scriptsize{res}}_G^H\chi_\nu}e(\eta)=\sum_{\eta\mid \mbox{\scriptsize{res}}_G^H\chi_\nu}\left(\frac{\eta(1)}{|H|}\sum_{h\in H}\eta(h^{-1})h\right)\in (\QQ^c\Gamma_0)[H],$$
 we achieve
$$e=\sum_\nu e_{\chi_\nu}\in A\cap (\QQ^c\Gamma_0)[H].$$
Since $e\in A\subseteq \QQ G$ is invariant under the action of $G(\Q_l^c/\Q_l)$, we get $e\in (\QQ\Gamma_0)[H]$ and finally
$$e\in A\cap(\QQ\Gamma_0)[H].$$
This proves the lemma.\hfill$\Box$

Now, $B:=A\cap(\QQ\Gamma_0)[H]$ is a two-sided ideal of $(\QQ\Gamma_0)[H]$ (because $A$ is a two-sided ideal in $\QQ G$) 
 and is thus the sum of some Wedderburn components of $(\QQ\Gamma_0)[H]$.  Because of $A e_\chi \neq 0$ and 
 $$(\QQ\Gamma_0)[H]e_\chi=(\QQ\Gamma_0)[H]\sum_{\eta\mid \mbox{\scriptsize{res}}_G^H\chi}e(\eta)\neq 0,$$
 there is a Wedderburn component $B_0$ of $(\QQ\Gamma_0)[H]$ in $B$ with $B_0 e_\chi \neq 0$. 
 Clifford theory (\cite[p.~565]{Hu}) shows that $$\res_G^H(\chi)=\sum_{j=0}^{w_\chi -1}\eta^{\gamma^{j}},$$
  where $\gamma^j$ runs through a set of representatives of $G/St(\eta)$.
  Thus, without loss of generality, we can assume $B_0 e(\eta)\neq 0$ because $e_\chi=\sum_{j=0}^{w_\chi-1}e(\eta^{\gamma^j})$. 
  
But $e(\eta)=\frac{\eta(1)}{|H|}\sum_{h\in H}\eta(h^{-1})h\notin (\QQ\Gamma_0)[H]$ because $\eta(H)\nsubseteq \Q_l$.  With $$\Gal:=G(\QQ^{\Q_l(\eta)}\Gamma_0/\QQ\Gamma_0),$$
 we therefore have 
$$B_0=(\QQ\Gamma_0)[H]\sum_{\sigma\in  \gal}e(\eta^\sigma)=:(\QQ\Gamma_0)[H]\varepsilon(\eta)$$
and $\varepsilon(\eta)$ is the central primitive idempotent of $B_0$.
Observe that, by the definition of $\QQ^{\Q_l(\eta)}\Gamma_0$, the Galois group $\Gal=G(\QQ^{\Q_l(\eta)}\Gamma_0/\QQ\Gamma_0)$ is canonically isomorphic to $G(\Q_l(\eta)/\Q_l)$. Thus, $\Gal$ is cyclic because $l\neq 2$. 

The other Wedderburn components of $B$ belong to the other irreducible characters $\eta^{\gamma^j}$ because these are exactly the characters not annihilating $B$. Thus, the structure of $B$ is given by
$$B=\bigoplus_{j=0}^{v_\chi-1}(\QQ\Gamma_0)[H]\varepsilon(\eta^{\gamma^j})$$
with 
$$v_\chi=\mbox{min}\{0\leq j\leq w_\chi-1 : \eta^{\gamma^j}=\eta^\sigma \mbox{ for some } \sigma\in \Gal\};$$
 note that $v_\chi\mid w_\chi$. We need $v_\chi$ to avoid that some $e(\eta^{\gamma^j})$ appears more than once as summand of the central primitive idempotents of $B$. Thus, the choice of $v_\chi$ ensures that the sum is direct.

By Roquette (\cite[Satz 2]{Ro}), the direct summands  of $B$ have trivial Schur index and are thus full matrix rings over certain fields. More precisely, we achieve
$$B\cong \bigoplus_{j=0}^{v_\chi-1}(\QQ^{\Q_l(\eta^{\gamma^j})}\Gamma_0)_{\eta^{\gamma^j}(1)\times\eta^{\gamma^j}(1)}=\bigoplus_{j=0}^{v_\chi-1}(\QQ^{\Q_l(\eta)}\Gamma_0)_{\eta(1)\times\eta(1)}$$
and
\begin{align*}
A&\stackrel{1}{=}\bigoplus_{i=0}^{l^m-1} B \gamma^{i}=\bigoplus_{i=0}^{l^m-1}\left(\bigoplus_{j=0}^{v_\chi-1}(\QQ\Gamma_0)[H]\varepsilon(\eta^{\gamma^j})\right)\gamma^{i}\\
&\cong \bigoplus_{i=0}^{l^m-1}\left(\bigoplus_{j=0}^{v_\chi-1}(\QQ^{\Q_l(\eta)}\Gamma_0)_{\eta(1)\times\eta(1)}\right)\gamma^{i}.
\end{align*}
Here, we have used for $\stackrel{1}{=}$ that, on the one hand, the direct sum is  contained in the $\QQ \Gamma_0$-algebra $A$  by the definition of $B$. On the other hand, this direct sum is, as a two-sided ideal of $\QQ G$, the direct sum of some Wedderburn components. Because $A$ is simple, the equality holds. 

Next, we define 
$$G_0:=\{\sigma\in \Gal:\eta^\sigma=\eta^{\gamma^j} \mbox{ for a } 0\leq j\leq w_\chi-1\}.$$
 This group is strongly related to $v_\chi$: The minimal $0\leq j\leq w_\chi-1$ satisfying the condition $\eta^\sigma=\eta^{\gamma^j}$ for a $\sigma \in \Gal$ is by definition $v_\chi$.

\begin{proposition}
Let $A$ be the simple component of $\QQ G$ corresponding to the irreducible character $\chi\in R_l G$. Then, 
$$A\cong \bigoplus_{i=0}^{l^m-1}\left(\bigoplus_{j=0}^{v_\chi-1}(\QQ^{\Q_l(\eta)}\Gamma_0)_{\eta(1)\times\eta(1)}\right)\gamma^{i}
$$
has centre 
$$Z(A)\cong \QQ^{L}\Gamma^{w_\chi}$$
with $L=\Q_l(\eta)^{G_0}$ and $G_0=\{\sigma\in \Gal:\eta^\sigma=\eta^{\gamma^j} \mbox{ for a } 0\leq j\leq w_\chi-1\}$, as before. \\
Moreover, $G_0$ is a cyclic group of order $\frac{w_\chi}{v_\chi}$.
\end{proposition}

\textbf{Proof:} 
First, we examine $G_0$. Define $\sigma_{v_\chi}\in \Gal$ via $\eta^{\gamma^{v_\chi}}=\eta^{\sigma_{v_\chi}}$. Observe that the actions of $\gamma$ and $\sigma\in \Gal$ commute. 
 Thus, induction yields $\eta^{\gamma^{n v_\chi}}=\eta^{\sigma_{v_\chi}^n}$ for $0\leq n\leq w_\chi/v_\chi-1$ and therefore, by $\eta^{\gamma^{v_\chi\cdot w_\chi/v_\chi}}=\eta^{\gamma^{w_\chi}}=\eta$, the order of $\sigma_{v_\chi}$ is $w_\chi/v_\chi$.
 
Furthermore, $G_0$ is, as subgroup of the cyclic group $\Gal\cong G(\Q_l(\eta)/\Q_l)$, cyclic. Next, every $0\leq j\leq w_\chi-1$ such that there exists a $\sigma\in \Gal$ with $\eta^\sigma=\eta^{\gamma^j}$ is a multiple of $v_\chi$, otherwise $v_\chi$ would not be minimal. Thus, we  conclude that $\sigma_{v_\chi}$ is the generator of $G_0$, we obtain $G_0=\langle \sigma_{v_\chi}\rangle$.
Next, we compute the centre of $A$. To do so, take a central element $z=\sum_{i=0}^{l^m-1} b_i\gamma^{i}$ with $b_i=\sum_{j=0}^{v_\chi-1}\beta _{ij}\varepsilon(\eta^{\gamma^j})\in B$
. Observe that the $i$-sum  is not orthogonal 
 and thus
$$\gamma z\stackrel{!}{=}z\gamma\Longleftrightarrow \sum_{i=0}^{l^m-1}\gamma b_i\gamma^{i}=\sum_{i=0}^{l^m-1} b_i^{\gamma^{-1}}\gamma^{i+1}\stackrel{!}{=}\sum_{i=0}^{l^m-1} b_i\gamma^{i+1}.$$
This is equivalent to $b_i^{\gamma^{-1}}=b_i$ for all $0\leq i\leq l^m-1$, i.e. 
$$\sum_{j=0}^{v_\chi-1}\beta _{ij}^{\gamma^{-1}}\varepsilon(\eta^{\gamma^{j-1}})=\sum_{j=0}^{v_\chi-1}\beta _{ij}\varepsilon(\eta^{\gamma^j})$$
 and thus, because we can read $j$ modulo $v_\chi$, we have
\begin{align}\label{1}
\beta_{i0}^{\gamma^{-1}}\varepsilon(\eta^{\gamma^{v_\chi-1}})=\beta_{i,v_\chi-1}\varepsilon(\eta^{\gamma^{v_\chi-1}}),\ ...\ ,\ \beta_{i1}^{\gamma^{-1}}\varepsilon(\eta)=\beta_{i0}\varepsilon(\eta).
\end{align}
This implies that  $\beta_{i0}$ determines $b_i$.

Analogously, $\gamma^{v_\chi}z\stackrel{!}{=}z\gamma^{v_\chi}$ and $\varepsilon(\eta^{\gamma^{v_\chi}})=\varepsilon(\eta)$ yield
\begin{align}\label{2}
\beta_{i0}^{\gamma^{-v_\chi}}\varepsilon(\eta)=\beta_{i0}\varepsilon(\eta)\ \forall\ 0\leq i\leq l^m-1.
\end{align}

Moreover, the central element $z$ commutes with  $\varepsilon(\eta)$: 
$$\varepsilon(\eta)z\stackrel{!}{=}z\varepsilon(\eta)\Longleftrightarrow\sum_{i=0}^{l^m-1} \varepsilon(\eta)b_i\gamma^{i}
\stackrel{!}{=}
\sum_{i=0}^{l^m-1} b_i\gamma^{i}\varepsilon(\eta)=
\sum_{i=0}^{l^m-1} b_i\varepsilon(\eta)^{\gamma^{-i}}\gamma^{i}.$$
This is equivalent to $\varepsilon(\eta)b_i=b_i\varepsilon(\eta)^{\gamma^{-i}}=b_i\varepsilon(\eta^{\gamma^{-i}})$  for all $0\leq i\leq l^m-1$. Because the primitive central idempotents $\varepsilon(\eta^{\gamma^{-i}})$ are orthogonal, i.e.~$\varepsilon(\eta^{\gamma^{-i}})\cdot \varepsilon(\eta^{\gamma^{-j}})\neq 0$ iff $i\equiv j\mod v_\chi$, we get $\varepsilon(\eta)b_i=0$
and thus 
$$b_i=0\ \forall\ v_\chi\nmid i$$  
by $\varepsilon(\eta)b_i=\varepsilon(\eta)\beta_{i0}$ and (\ref{1}).
From now on, we therefore consider  $z=\sum'b_i\gamma^{i}$ with ${\sum}':=\sum_{i=0,v_\chi\mid i}^{l^m-1}$ .

Finally, the central element $z$ commutes with every $b\varepsilon(\eta)$ for $b\in (\QQ\Gamma_0)[H]$:
$$b\varepsilon(\eta)z\stackrel{!}{=}zb\varepsilon(\eta)\Longleftrightarrow {\sum}' b\varepsilon(\eta)b_i\gamma^{i}
\stackrel{!}{=}{\sum}' b_i\gamma^{i}b\varepsilon(\eta)=
{\sum}' b_ib^{\gamma^{-i}}\varepsilon(\eta)\gamma^{i}$$
and thus 
\begin{align}\label{3}
b\varepsilon(\eta)b_i&=(b\varepsilon(\eta))(b_i\varepsilon(\eta))=(b\varepsilon(\eta))(\beta_{i0}\varepsilon(\eta))=(\beta_{i0}\varepsilon(\eta))(b^{\gamma^{-i}}\varepsilon(\eta))
\end{align}
for all $b\in (\QQ\Gamma_0)[H]$ and for all $0\leq i\leq l^m-1,\ v_\chi\mid i$.

To understand this condition, we first consider $w_\chi\mid i$. In this case, the representation 
$$D:(\QQ\Gamma_0)[H]\rightarrow (\QQ^{\Q_l(\eta)}\Gamma_0)_{\eta(1)\times\eta(1)}$$
 corresponding to $\eta$ leads to
$$D(b)D(\beta_{i0})=D(\beta_{i0})D(b^{\gamma^{-i}}).$$
 From the definition of $w_\chi$, we achieve $D(b^{\gamma^{-w_\chi}})=YD(b)Y^{-1}$ and $$D(b^{\gamma^{-i}})=Y^{i/w_\chi}D(b)Y^{-(i/w_\chi)}=:Y_iD(b)Y_i^{-1}$$ with a matrix $Y\in (\QQ^{\Q_l(\eta)}\Gamma_0)_{\eta(1)\times\eta(1)}$. Observe that $Y_i$ is unique up to central elements of $(\QQ^{\Q_l(\eta)}\Gamma_0)_{\eta(1)\times\eta(1)}$. Moreover, $Y$ is independent of $b$ because $D$ is multiplicative and thus, even though $\gamma^{w_\chi}\notin (\QQ \Gamma_0)[H]$, we see $$D((bc)^{\gamma^{-w_\chi}})=D(\gamma^{w_\chi}b\gamma^{-w_\chi}\gamma^{w_\chi}c\gamma^{-w_\chi})=D(b^{\gamma^{-w_\chi}})D(c^{\gamma^{-w_\chi}}).$$
  This yields
\begin{align*}
D(b)D(\beta_{i0})=D(\beta_{i0})Y_iD(b)Y_i^{-1}\ \Longleftrightarrow\ D(b)(D(\beta_{i0})Y_i)=(D(\beta_{i0})Y_i)D(b).
\end{align*}
Since $D$ is surjective, $D(b)$ runs through $(\QQ^{\Q_l(\eta)}\Gamma_0)_{\eta(1)\times\eta(1)}$ and $D(\beta_{i0})Y_i$ is therefore   central in $(\QQ^{\Q_l(\eta)}\Gamma_0)_{\eta(1)\times\eta(1)}$.   Furthermore, there exists a $y_i\in (\QQ\Gamma_0)[H]$ with $D(y_i)=Y_i$ and $y_i$ can be chosen as $y_i=y_{w_\chi}^{i/w_\chi}$ because $D$ is an epimorphism. 
Thus, for every $0\leq i\leq l^m-1$ with $w_\chi\mid i$, 
the element $(\beta_{i0}\varepsilon(\eta))(y_i\varepsilon(\eta))$  maps to a central matrix in $\QQ^{\Q_l(\eta)}\Gamma_0\cdot \mbox{\textbf{1}}$. 

Because of $\gamma^{v_\chi}\varepsilon(\eta)\in A$, the centre has to be fixed under the action induced by the conjugation by $\gamma^{v_\chi}$, which by construction is equal to the action of  $\sigma_{v_\chi}\in G_0=\langle\sigma_{v_\chi}\rangle$ on $\QQ^{\Q_l(\eta)}\Gamma_0$, i.e.
\begin{align}\label{a}
(\beta_{i0}\varepsilon(\eta))(y_i\varepsilon(\eta))\mapsto D(\beta_{i0})D(y_i) \in \QQ^L\Gamma_0\cdot \mbox{\textbf{1}}.
\end{align} 

Next, we show that $b_i=0$ for $0\leq i\leq l^m-1$ with $v_\chi\mid i$ but $w_\chi\nmid i$. Observe that conjugation by $\gamma^{v_\chi}$ induces an automorphism on  $(\QQ\Gamma_0)[H]\varepsilon(\eta)\cong (\QQ^{\Q_l(\eta)}\Gamma_0)_{\eta(1)\times\eta(1)}$ because $H$ is normal in $G$ and $\varepsilon(\eta)^{\gamma^{v_\chi}}=\varepsilon(\eta^{\gamma^{v_\chi}})=\varepsilon(\eta)$. This automorphism fixes the centre  of $(\QQ\Gamma_0)[H]\varepsilon(\eta)$ because for $z$  a central element and $x\in (\QQ\Gamma_0)[H]\varepsilon(\eta)$, we have
$$\gamma^{-v_\chi}z\gamma^{v_\chi}x=\gamma^{-v_\chi}z\gamma^{v_\chi}x\gamma^{-v_\chi}\gamma^{v_\chi}=\gamma^{-v_\chi}\gamma^{v_\chi}x\gamma^{-v_\chi}z\gamma^{v_\chi}=x\gamma^{-v_\chi}z\gamma^{v_\chi}.$$
Thus, it fixes the centre of $(\QQ^{\Q_l(\eta)}\Gamma_0)_{\eta(1)\times\eta(1)}$, too. But this automorphism does not fix\\ $(\QQ\Gamma_0)[H]\varepsilon(\eta)\cong(\QQ^{\Q_l(\eta)}\Gamma_0)_{\eta(1)\times\eta(1)}$ elementwise; we will call it $c_{\gamma^{v_\chi}}$.
 $\sigma_{v_\chi}$ can also be extended to an automorphism on $(\QQ^{\Q_l(\eta)}\Gamma_0)_{\eta(1)\times\eta(1)}$ via its action on every entry of the matrices in $(\QQ^{\Q_l(\eta)}\Gamma_0)_{\eta(1)\times\eta(1)}$. This yields the central automorphism $c_{\gamma^{v_\chi}}\sigma_{v_\chi}^{-1}$ of $ (\QQ^{\Q_l(\eta)}\Gamma_0)_{\eta(1)\times\eta(1)}$ which,  by the theorem of Skolem-Noether, is 
 the conjugation by a matrix $\tilde{Y}\in (\QQ^{\Q_l(\eta)}\Gamma_0)_{\eta(1)\times\eta(1)}$. 

Now, we can repeat the above calculations for this case. Let $\sigma_{i}$ be the automorphism induced by $\gamma^{i}$ and observe that $\sigma_i=\sigma_{v_\chi}^{i/v_\chi}$.
This time, we get 
\begin{align}\label{111}
D(b^{\gamma^{-i}})^{\sigma_i}=\tilde{Y}_iD(b)\tilde{Y}_i^{-1} \ \Longleftrightarrow\  D(b^{\gamma^{-i}})=\tilde{Y}_i^{\sigma_i^{-1}}D(b)^{\sigma_i^{-1}}\tilde{Y}_i^{\sigma_i^{-1}-1}.
\end{align}
With $Y_i:=\tilde{Y}_i^{\sigma_i^{-1}}$, condition (\ref{3}) yields
\begin{align}\label{5}
D(b)D(\beta_{i0})=D(\beta_{i0})Y_iD(b)^{\sigma_i^{-1}}Y_i^{-1}\ \Longleftrightarrow\ 
D(b)(D(\beta_{i0})Y_i)=(D(\beta_{i0})Y_i)D(b)^{\sigma_i^{-1}}.
\end{align}
Because $b\in (\QQ\Gamma_0)[H]$ is arbitrary and $D$ is surjective, we can choose  $b\in (\QQ\Gamma_0)[H]$ such that $D(b)=\alpha\cdot\mbox{\bfseries{1}}$ is central and $\alpha^{\sigma_i^{-1}}\neq \alpha\in \QQ^{\Q_l(\eta)}\Gamma_0$. Now, (\ref{5}) yields $$\alpha\mbox{\bfseries{1}}\cdot D(\beta_{i0})Y_i= D(\beta_{i0})Y_i \alpha^{\sigma_i^{-1}}\mbox{\bfseries{1}}=\alpha^{\sigma_i^{-1}}\mbox{\bfseries{1}}\cdot D(\beta_{i0})Y_i$$
 and thus $\alpha\mbox{\bfseries{1}}\cdot D(\beta_{i0})=\alpha^{\sigma_i^{-1}}\mbox{\bfseries{1}}\cdot D(\beta_{i0})$ by cancelling $Y_i$. Therefore, we conclude $\alpha x_{\nu\tau}=\alpha^{\sigma_i^{-1}}x_{\nu\tau}$ for every entry $x_{\nu\tau}$ of $D(\beta_{i0})$ and finally $D(\beta_{i0})=0$.

These results are now summarized. Let $z=\sum' b_i\gamma^{i}$ with $b_i=\sum_{j=0}^{v_\chi-1}\beta _{ij}\varepsilon(\eta^{\gamma^j})\in B$ be a central element of $A$. We have seen that we can assume $\sum'=\sum_{i=0,w_\chi\mid i}^{l^m-1}$ and that $\beta_{i0}$ determines $b_i$ uniquely for each such $i$. 

In remains to show the claimed isomorphism for $Z(A)$. We start with the proof that
\begin{align*}
\varphi_1: Z(A)&\rightarrow {\bigoplus}'(\QQ\Gamma_0)[H]\varepsilon(\eta)\gamma^{i},\\
 {\sum}' b_i\gamma^{i}&\mapsto {\sum}'b_i\varepsilon(\eta)\gamma^{i}={\sum}'\beta_{i0}\varepsilon(\eta)\gamma^{i},
\end{align*}
is a monomorphism of $\QQ\Gamma_0$-algebras. First, $\varphi_1$ is the identity on $\QQ\Gamma_0$. Next, the map $\varphi_1$ ist injective because $\beta_{i0}\varepsilon(\eta)$ determines $b_i$ by condition (\ref{1}). Additivity of $\varphi_1$ being obvious, it remains to check  multiplicity: Let $ \sum' b_i\gamma^{i}$ and $\sum' \tilde{b}_j\gamma^{j}\in Z(A)$. Then
\begin{align*}
\varphi_1\left({\sum}' b_i\gamma^{i}\right)\cdot\varphi_1\left({\sum}' \tilde{b}_j\gamma^{j}\right)
&=\left({\sum}'\beta_{i0}\varepsilon(\eta)\gamma^{i}\right)\cdot\left({\sum}'\tilde{\beta}_{j0}\varepsilon(\eta)\gamma^{j}\right)\\
&\stackrel{1}{=}{\sum}'\left(\sum_{\nu+\tau=i}\ '\beta_{\nu 0}(\tilde{\beta_{\tau 0}})^{\gamma^{-\nu}}\right)\varepsilon(\eta)\gamma^{i}.\\
&\stackrel{2}{=}{\sum}'\left({\sum_{\nu+\tau=i}}'\beta_{\nu 0}\tilde{\beta_{\tau 0}}\right)\varepsilon(\eta)\gamma^{i}.
\end{align*}
For $\stackrel{1}{=}$, we have used the fact that $\varepsilon(\eta)$ is a central idempotent of $(\QQ\Gamma_0)[H]$ and $\varepsilon(\eta^{w_\chi})=\varepsilon(\eta)$; for $\stackrel{2}{=}$ use condition (\ref{2}). Finally 
\begin{align*}
\varphi_1\left(\left({\sum}' b_i\gamma^{i}\right)\cdot\left({\sum}' \tilde{b}_j\gamma^{j}\right)\right)
&=\varphi_1\left({\sum}' \left({\sum_{\nu+\tau=i}}'b_\nu \tilde{b}_\tau\right)\gamma^{i}\right)\\
&={\sum}' \left({\sum_{\nu+\tau=i}}'b_\nu \tilde{b}_\tau\right)\varepsilon(\eta)\gamma^{i}\\
&={\sum}' \left({\sum_{\nu+\tau=i}}'\beta_{\nu0} \tilde{\beta}_{\tau0}\right)\varepsilon(\eta)\gamma^{i}
\end{align*}
and multiplicity is stated.

The action of $\gamma^{w_\chi}=(\gamma^{v_\chi})^{w_\chi/v_\chi}$ on $L=\Q_l(\eta)^{G_0}$ is trivial by the definition of $v_\chi$. It is clearly trivial on $\QQ\Gamma_0$, too. This yields the identity
$${\bigoplus}'(\QQ^{L}\Gamma_0)\gamma^{i}=\QQ^{L}\Gamma^{w_\chi}.$$
 Next, we fix $y_{w_\chi}\in (\QQ \Gamma_0)[H]$ and $y_i=y_{w_\chi}^{i/w_\chi}$ as above and claim that 
\begin{align*}
\varphi: Z(A)&\rightarrow {\bigoplus}'(\QQ^{L}\Gamma_0)\gamma^{i}\cdot \mbox{\textbf{1}}=\QQ^{L}\Gamma^{w_\chi}\cdot \mbox{\textbf{1}},\\ {\sum}' b_i\gamma^{i}&\mapsto {\sum}'D(\beta_{i0})D(y_i)\gamma^{i},
\end{align*}
 is an isomorphism of $\QQ\Gamma_0$-algebras (compare (\ref{a})). Again, $\varphi$ is the identity on $\QQ\Gamma_0$. We can write $\varphi$ as 
$$\varphi=D\circ r_{y} \circ \varphi_1$$
with 
\begin{align*}
r_{y}:{\bigoplus}'(\QQ\Gamma_0)[H]\varepsilon(\eta)\gamma^{i}&\rightarrow {\bigoplus}'(\QQ\Gamma_0)[H]\varepsilon(\eta)\gamma^{i},\\
 {\sum}' \beta_{i0}\varepsilon(\eta)\gamma^{i}&\mapsto {\sum}' \beta_{i0}\varepsilon(\eta) y_i\varepsilon(\eta)\gamma^{i},
\end{align*}
a monomorphism of additive groups.

The representation $D$ is here extended to ${\bigoplus}'(\QQ\Gamma_0)[H]\varepsilon(\eta)\gamma^{i}$ by $$D\left({\sum}'\beta_{i0}\varepsilon(\eta)\gamma^{i}\right)={\sum}'D(\beta_{i0})\gamma^{i},$$
 which again is a homomorphism. Observe that, in fact, we only need to extend $D|_{(\QQ \Gamma_0)[H]\varepsilon(\eta)}$. Thus, the additivity of $\varphi$ results from the additivity of $\varphi_1$, $r_{y}$ and $D$, the injectivity of $\varphi$ follows from the injectivity of $\varphi_1$, $r_y$ and $D|_{(\QQ\Gamma_0)[H]\varepsilon(\eta)}$. 

For multiplicity, we have on the one hand
\begin{align*}
\varphi\left({\sum}' b_i\gamma^{i}\right)\cdot\varphi\left({\sum}'\tilde{b}_i\gamma^{i}\right)
&=\left({\sum}'D(\beta_{i0})D(y_i)\gamma^{i}\right)\cdot\left({\sum}' D(\tilde{\beta}_{i0})D(y_i)\gamma^{i}\right)\\
&\stackrel{1}{=}{\sum}' \left({\sum_{\nu+\tau=i}}'D(\beta_{\nu0})D(y_\nu) D(\tilde{\beta}_{\tau0})D(y_\tau)\right)\gamma^{i}\\
&\stackrel{2}{=}{\sum}' \left({\sum_{\nu+\tau=i}}'D(\beta_{\nu0}) D(\tilde{\beta}_{\tau0})D(y_\nu)D(y_\tau)\right)\gamma^{i}\\
&\stackrel{3}{=}{\sum}' \left({\sum_{\nu+\tau=i}}'D(\beta_{\nu0}) D(\tilde{\beta}_{\tau0})\right)D(y_i)\gamma^{i}.\\
\end{align*}
Here, $\stackrel{1}{=}$ is condition (\ref{a}) together with the fact that $w_\chi\mid i$. For $\stackrel{2}{=}$,  we have used that $D(y_\nu)=Y_\nu$ commutes with $D(\beta_{\tau0})$ by condition (\ref{2}) and the definition of $Y_\nu$.
 For $\stackrel{3}{=}$, observe that $$D(y_\nu)D(y_\tau)=Y_{w_\chi}^{(\nu+\tau)/w_\chi}=D(y_{\nu+\tau})=D(y_i).$$
 
 On the other hand, we keep in mind that $\tilde{b}_i$ and $\gamma^i$ commute by (\ref{2}) and compute
\begin{align*}
\varphi\left(\left({\sum}' b_i\gamma^{i}\right)\cdot\left({\sum}' \tilde{b}_i\gamma^{i}\right)\right)
&=\varphi\left({\sum}' \left({\sum_{\nu+\tau=i}}'b_\nu \tilde{b}_\tau\right)\gamma^{i}\right)\\
&=(D\circ r_{y})\left({\sum}' \left({\sum_{\nu+\tau=i}}'\beta_{\nu0}\tilde{\beta}_{\tau 0}\varepsilon(\eta) \right)\gamma^{i}\right)\\
&=D\left({\sum}' \left({\sum_{\nu+\tau=i}}'\beta_{\nu0}\tilde{\beta}_{\tau 0}\varepsilon(\eta)  \right)y_i\varepsilon(\eta)\gamma^{i}\right)\\
&={\sum}' D\left({\sum_{\nu+\tau=i}}'\beta_{\nu0} \tilde{\beta}_{\tau0}\right)D(y_i)\gamma^{i}\\
&={\sum}' \left({\sum_{\nu+\tau=i}}'D(\beta_{\nu0}) D(\tilde{\beta}_{\tau0})\right)D(y_i)\gamma^{i}.
\end{align*}
Thus, multiplicity is stated.

It remains to prove that $\varphi$ is surjective. 
To do so, take a $W_i\in \QQ^L\Gamma_0\cdot\mbox{\textbf{1}}$ for every $0\leq i\leq l^m-1$ with $w_\chi\mid i$. By the surjectivity of $D$, there exists a $w_i\in (\QQ \Gamma_0)[H]\varepsilon(\eta)$ such that $W_i=D(w_i)$. Define $\beta_{i0}\in (\QQ\Gamma_0)[H]$ via $D(w_i)=D(\beta_{i0})D(y_i)$. This $\beta_{i0}$ determines $b_i$ via (\ref{1}) uniquely and finally $z=\sum{'}b_i\gamma^{i}$ maps to $\sum{'}W_i\gamma^{i}$.\hfill$\Box$

Now, we can compute the structure of $A$.

\begin{theorem}\label{th1}
Let $A$ be the simple component of $\QQ G$ corresponding to the irreducible $\chi\in R_lG$ and $\sigma_{v_\chi}\in G_0$ with $\langle\sigma_{v_\chi}\rangle=G_0$ and $\eta^{\sigma_{v_\chi}}=\eta^{\gamma^{v_\chi}}$.
Then 
\begin{enumerate}
\item $\dim_{Z(A)}A=\chi(1)^2$,
\item $A$ is split by $\QQ^{\Q_l(\eta)}\Gamma^{w_\chi}$, 
\item $A$ has Schur index $s_D=\frac{w_\chi}{v_\chi}$ and 
\item $A\cong D_{n\times n}$ with $n=\frac{\chi(1)}{s_D}$ and the skew field $D$ is cyclic:
$$D\cong  \bigoplus_{i=0}^{w_\chi/v_\chi-1} (\QQ^{\Q_l(\eta)}\Gamma^{w_\chi}) \gamma^{v_\chi i}=: (\QQ^{\Q_l(\eta)}\Gamma^{w_\chi}/\QQ^L\Gamma^{w_\chi},\sigma_{v_\chi},\gamma^{w_\chi})$$
with $L=\Q_l(\eta)^{G_0}$.
\end{enumerate}
\end{theorem}

\textbf{Proof:}
For (i), let $\chi'$ be an irreducible $\Q_l^c$-character of $St(\eta)$ extending $\eta$ such that $\chi=\ind_{St(\eta)}^G(\chi')$.
  Consider $Z_{\chi'}(A):=Z(A)(\chi')$. From
\begin{align}\label{8}
\Q_l(\chi')\otimes_{\Q_l}A=(\Q_l(\chi')\otimes_{\Q_l}Z(A))\otimes_{Z(A)}A\cong  (Z_{\chi'}(A)\otimes_{Z(A)}A)^{r_{\chi'}}
\end{align}
with $r_{\chi'}$ appropriately, we conclude $Z_{\chi'}(A)\cong Z(\QQ^{\Q_l(\chi')}Ge_\chi)$: 
First, we show
\begin{align}\label{7}
\Q_l(\chi')\otimes_{\Q_l}A={\bigoplus}_{\sigma\in \gal(\chi')}^* (\QQ^{\Q_l(\chi')}G)e_{\chi^\sigma}
\end{align}
with $\Gal(\chi'):=G(\QQ^{\Q_l(\chi')}\Gamma_0/\QQ\Gamma_0)$ and $\bigoplus^*$ meaning summation modulo type $W$.

Because $\chi$ can be regarded as irreducible character of the finite group $G/\Gamma_0$, there exists a representation of $\chi$ over the fields $\Q_l(\chi)\subseteq \Q_l(\chi')$ by \cite[Satz 1]{Ro}. 
Furthermore, the irreducible constituents $\eta_j$ of $\res_G^H\chi^\sigma$ are precisely the characters $(\eta^{\gamma^j})^\sigma$ obtained by the irreducible constituents $\eta^{\gamma^j}$ of $\res_G^H\chi$. We further recall that the actions of $\gamma$ and $\sigma$ commute. Thus, $$e_{\chi^\sigma}=\sum_{\eta_j\mid \scriptsize{\res}_G^H\chi^\sigma}e(\eta_j)=\sum_{j=0}^{w_\chi-1}e(\eta^{\gamma^j})^{\sigma}=e_\chi^\sigma \in \QQ^{\Q_l(\chi')}G$$
 because $\Q_l(\eta)\subseteq \Q_l(\chi')$ by construction and $e_{\chi^\sigma}$ is therefore a primitive central idempotent of $\QQ^{\Q_l(\chi')}G$, i.e.~
 $\QQ^{\Q_l(\chi')}Ge_{\chi^\sigma}$ is a simple component of $\QQ^{\Q_l(\chi')}G$.  
 Hence summation modulo type W ensures that the summands of the right hand side are distinct because in \cite[Cor, p. 556]{Iw2} it is shown that two primitive central idempotents are equal iff the corresponding irreducible characters only differ by a $W$-twist.

The inclusion `$\supseteq$' of (\ref{7}) is true by the following.
From \cite[Satz 1]{Ro}, we know that
 $\eta$ has a realization over $\QQ^{\Q_l(\eta)}\Gamma_0$ and thus over $\QQ^{\Q_l(\chi')}\Gamma_0\supseteq \QQ^{\Q_l(\eta)}\Gamma_0$.     This yields  that the $(\eta^{\gamma^j})^\sigma$ are absolute irreducible characters of $H$ belonging to  Wedderburn components of $\Q_l(\chi')\otimes_{\Q_l}B$, recall that $B=\bigoplus_{j=0}^{v_\chi-1}(\QQ\Gamma_0)[H]\varepsilon(\eta^{\gamma^j})$. Therefore,
 $Be_{\chi^\sigma}\neq 0$ and in particular 
  $Ae_{\chi^\sigma}\neq 0$. Finally, $\QQ^{\Q_l(\chi')}Ge_{\chi^\sigma}\subseteq \Q_l(\chi')\otimes_{\Q_l}A$ follows. As the Wedderburn components are orthogonal, the direct sum of the $\QQ^{\Q_l(\chi')}Ge_{\chi^\sigma}$ is also contained in $\Q_l(\chi')\otimes_{\Q_l}A$ and the claimed inclusion follows.

For the other inclusion,  observe that $\Q_l(\chi')\otimes_{\Q_l}A\subseteq \QQ^{\Q_l(\chi')}G$ carries a natural $\Gal(\chi')$-action with fixed points $A\subset \QQ G$. Also, $\bigoplus_{\sigma\in \gal(\chi')}^*(\QQ^{\Q_l(\chi')}G)e_{\chi^\sigma}$ carries this action. As the set of fixed points of $\QQ^{\Q_l(\chi')}G=\Q_l(\chi')\otimes_{\Q_l} \QQ G$ is exactly $\QQ G$, the fixed points of the right side make up a two-sided ideal of $\QQ G$. 
Of course, the fixed points of $\bigoplus_{\sigma\in \gal(\chi')}^*(\QQ^{\Q_l(\chi')}G)e_{\chi^\sigma}\subseteq \Q_l(\chi')\otimes_{\Q_l}A$ are contained in the set of fixed points of $\Q_l(\chi')\otimes_{\Q_l}A$, i.e.~in $A$.
Since $A$ is simple, the fixed points in $\bigoplus_{\sigma\in \gal(\chi')}^*(\QQ^{\Q_l(\chi')}G)e_{\chi^\sigma}$ are therefore given by $A$. 
This implies the equality.

Thus, the simple components on the left and right side of (\ref{7}) coincide. Observe that $Z_{\chi'}(A)\otimes_{Z(A)} A$ is a central simple $Z_{\chi'}(A)$-algebra because $A$ is a central simple $Z(A)$-algebra (compare \cite[Hilfssatz 14.2]{Hu}). Together with (\ref{8}), this yields 
$$Z_{\chi'}(A)\otimes_{Z(A)} A\cong \QQ^{\Q_l(\chi')}G e_{\chi^\sigma}.$$ Now, we compute the centres of the simple components of (\ref{7}). On the right side, the centres are obviously the $Z(\QQ^{\Q_l(\chi')}G e_{\chi^\sigma})$. The simple components of the left side are, by (\ref{8}), isomorphic to the central simple $Z_{\chi'}(A)$-algebra $Z_{\chi'}(A)\otimes_{Z(A)} A$. Thus, (\ref{7}) implies the claimed isomorphism 
$$Z_{\chi'}(A)\cong Z(\QQ^{\Q_l(\chi')}G e_{\chi^\sigma}).$$

In \cite[p.~555]{Iw2}, it is shown that  $\dim_{Z(\QQ^{\Q_l(\chi')}Ge_{\chi})}\QQ^{\Q_l(\chi')}Ge_{\chi}=\chi(1)^2$
and thus
\begin{align*}
\dim_{Z(A)}A=\dim_{Z_{\chi'}(A)}(Z_{\chi'}(A)\otimes_{Z(A)}A)=\chi(1)^2
\end{align*}
results.

We turn to the proof of (iv). The theorem of Wedderburn  implies that $A\cong D_{n\times n}$ is a full matrix ring over a skew field $D$. The dimension of $D$ over the centre of $A$ is the square $s_D^2$ with $s_D$ being the Schur index of $A$. Thus 
$$\chi(1)^2=\dim_{Z(A)}A=s_D^2n^2.$$
For the computation of $D$, we use the fact that $A\cong E^n$ for a minimal right ideal $E=\varepsilon A$ of $A$ and $D\cong \varepsilon A\varepsilon$ with $\varepsilon$ a primitive idempotent of $A$.
Analogously, there exists a primitive idempotent $\varepsilon_1$ of $(\QQ\Gamma_0)[H]\varepsilon(\eta)$ for the minimal right ideal $S_{\eta}=\varepsilon_1(\QQ\Gamma_0)[H]\varepsilon(\eta)$ of  $(\QQ\Gamma_0)[H]\varepsilon(\eta)$. We have seen that $(\QQ\Gamma_0)[H]\varepsilon(\eta)$ is the full matrix ring $(\QQ^{\Q_l(\eta)}\Gamma_0)_{\eta(1)\times\eta(1)}$
 and therefore $\varepsilon_1 (\QQ\Gamma_0)[H]\varepsilon(\eta)\varepsilon_1\cong \QQ^{\Q_l(\eta)}\Gamma_0$ with 
\[\varepsilon_1\mapsto\left(
\begin{array}{cccc}
1&0&\ldots&0\\
0&0&\ldots&0\\
\vdots&\vdots&&\vdots\\
0&0&\ldots&0\end{array}
\right)\in (\QQ^{\Q_l(\eta)}\Gamma_0)_{\eta(1)\times \eta(1)}.\]
 With $R:=(\QQ^{\Q_l(\eta)}\Gamma_0)_{\eta(1)\times\eta(1)}\cong(\QQ\Gamma_0)[H]\varepsilon(\eta)$, 
 we get
\begin{align*}
\bigoplus_{i=0}^{l^m-1}(S_\eta)\gamma^{i}
&=\bigoplus_{i=0}^{l^m-1}(\varepsilon_1R)\gamma^{i}
=\bigoplus_{i=0}^{l^m-1}\left(\underbrace{(\varepsilon_1,0,...,0)}_{v_\chi}\bigoplus_{j=0}^{v_\chi-1}R\right)\gamma^{i}\\
&=(\varepsilon_1,0,...,0)\bigoplus_{i=0}^{l^m-1}\left(\bigoplus_{j=0}^{v_\chi-1}R\right)\gamma^{i}\\
&\cong(\varepsilon_1,0,...,0)\bigoplus_{i=0}^{l^m-1}\left(\bigoplus_{j=0}^{v_\chi-1}(\QQ \Gamma_0)[H]\varepsilon(\eta^{\gamma^j})\right)\gamma^{i}
\end{align*}
and furthermore
\begin{align*}
\varepsilon_1A\varepsilon_1&\cong(\varepsilon_1,0,...,0)\left(\bigoplus_{i=0}^{l^m-1}\left(\bigoplus_{j=0}^{v_\chi-1}R\right)\gamma^{i}\right)(\varepsilon_1,0,...,0)\\
&=\bigoplus_{i=0}^{l^m-1}\left((\varepsilon_1,0,...,0)\bigoplus_{j=0}^{v_\chi-1}R(\varepsilon_1^{\gamma^{-i}},0,...,0)\right)\gamma^{i}\\
&=\bigoplus_{i=0}^{l^m-1}\left(\varepsilon_1R\varepsilon_1^{\gamma^{-i}}\right)\gamma^{i}\stackrel{1}{=}\bigoplus_{i=0, v_\chi\mid i}^{l^m-1}\left(\varepsilon_1R\varepsilon_1^{\gamma^{- i}}\right)\gamma^{i}.
\end{align*}
For $\stackrel{1}{=}$, we have used that conjugation by $\gamma^{i}$ permutes the Wedderburn components of $(\QQ\Gamma_0)[H]$, i.e.~$\varepsilon_1^{\gamma^{-i}}\in (Q\Gamma_0)[H]\varepsilon(\eta^{\gamma^{-i}})=(Q\Gamma_0)[H]\varepsilon(\eta^{\gamma^{v_\chi-i}})$, and thus only the conjugations by $\gamma^{i}$, $v_\chi\mid i$, yield nonzero summands. We have seen in (\ref{111}) that 
$$D(\varepsilon_1^{\gamma^{-i}})=Y_i D(\varepsilon_1)^{\sigma_i^{-1}} Y_i^{-1}$$
Therefore, $D(\varepsilon_1^{\gamma^{-i}})$ (resp.~$\varepsilon_1^{\gamma^{-i}}$) is still an idempotent of a minimal right ideal of $(\QQ^{\Q_l(\eta)}\Gamma_0)_{\eta(1)\times\eta(1)}\cong(\QQ\Gamma_0)[H]\varepsilon(\eta)$. Furthermore, it is a primitive idempotent because otherwise $(Y_i^{-1}D(\varepsilon_1^{\gamma^{-i}})Y_i)^{\sigma_i}= D(\varepsilon_1)$ is not primitive, a contradiction.

 As primitive idempotent, $\varepsilon_1^{\gamma^{-i}}$ is a matrix of rank $1$ and thus 
$$\varepsilon_1R\varepsilon_1^{\gamma^{-i}}
=\left\{\left(\begin{array}[]{ccc}
a_1&...&a_{\eta(1)}\\
0&...&0\\
\vdots&&\vdots\\0&...&0	
\end{array}\right)\alpha: \alpha\in \QQ^{\Q_l(\eta)}\Gamma_0\right\}
$$
with at least one of the $a_\nu\neq 0$ for $1\leq \nu\leq \eta(1)$. This implies  $\varepsilon_1R\varepsilon_1^{\gamma^{-i}}\cong \QQ^{\Q_l(\eta)}\Gamma_0$ as additive groups and also $$\varepsilon_1 A\varepsilon_1\cong \bigoplus_{i=0,v_\chi\mid i}^{l^m-1}\QQ^{\Q_l(\eta)}\Gamma_0 \gamma^{i}=\bigoplus_{i=0,v_\chi\mid i}^{w_\chi-1}\QQ^{\Q_l(\eta)}\Gamma^{w_\chi} \gamma^{i}=:A_0$$
as additive groups.
Concerning multiplication in $\varepsilon_1 A\varepsilon_1$, observe
\begin{align*}
\varepsilon_1R\varepsilon_1^{\gamma^{-i}}\gamma^{i}\cdot \varepsilon_1R\varepsilon_1^{\gamma^{-i'}}\gamma^{i'}
&\stackrel{1}{=}\varepsilon_1R\varepsilon_1^{\gamma^{-i}}\gamma^i R\varepsilon_1^{\gamma^{-i'}}\gamma^{i'}\\
&\stackrel{2}{=}\varepsilon_1R\varepsilon_1^{\gamma^{-i}} R\varepsilon_1^{\gamma^{-(i+i')}}\gamma^{i+i'}
\stackrel{3}{=}\varepsilon_1R\varepsilon_1^{\gamma^{-(i+i')}}\gamma^{i+i'}
\end{align*}
and $\stackrel{1}{=}$ is true because $\varepsilon_1^{\gamma^{-i}}$ is an idempotent, $\stackrel{2}{=}$ follows by $R\cong(\QQ \Gamma_0)[H]\varepsilon(\eta)$: As $H$ is normal in $G$, conjugation by $\gamma^{(-i)}$ fixes the algebra $(\QQ \Gamma_0)[H]$. Furthermore, $v_\chi\mid i$ ensures that even the Wedderburn components of $(\QQ \Gamma_0)[H]$ are fixed, i.e.~$R^{\gamma^{-i}}=R$. For $\stackrel{3}{=}$ note that $R$ is simple; thus, the two-sided ideal $R\varepsilon^{\gamma^{-i}}R$ is the whole ring $R$.   

 The direct $i$-sum in the upper description of $\varepsilon_1A\varepsilon_1$ is hence not orthogonal and the multiplication rules on $\varepsilon_1A\varepsilon_1$ and on $A_0$ coincide. Therefore, $\varepsilon_1A\varepsilon_1\cong A_0$ as $\QQ\Gamma_0$-algebras and $\varepsilon_1 A\varepsilon_1$ is the claimed crossed product $$\varepsilon_1A\varepsilon_1\cong(\QQ^{\Q_l(\eta)}\Gamma^{w_\chi}/\QQ^L\Gamma^{w_\chi},\sigma_{v_\chi},\gamma^{w_\chi}).$$

For $\varepsilon_1A\varepsilon_1$ to be a skew field, it remains to prove that the Schur index $s$ of the crossed product is equal to $[\QQ^{\Q_l(\eta)}\Gamma^{w_\chi}:\QQ^L\Gamma^{w_\chi}]=[\Q_l(\eta):L]=w_\chi/v_\chi$.  In this case, $\varepsilon_1A \varepsilon_1\cong D$ is the skew field underlying $A$ and $s_D=s=w_\chi/v_\chi$.

To show this, let $\NN$ denote the norm of the field extension $\QQ^{\Q_l(\eta)}\Gamma^{w_\chi}/\QQ^{L}\Gamma^{w_\chi}$ and set $o(\gamma^{w_\chi})$ the order of $\gamma^{w_\chi}$ in $(\QQ^{L}\Gamma^{w_\chi})^\times/\NN((\QQ^{\Q_l(\eta)}\Gamma^{w_\chi})^\times)$. We use the fact that  the crossed product $A_0$ is a skew field if $o(\gamma^{w_\chi})$ equals the degree $[\QQ^{\Q_l(\eta)} \Gamma^{w_\chi}:\QQ^L\Gamma^{w_\chi}]$ (see \cite[(30.7)]{Re}). Then, the Schur index is $s=[\QQ^{\Q_l(\eta)} \Gamma^{w_\chi}:\QQ^L\Gamma^{w_\chi}]=w_\chi/v_\chi$ as claimed.

For an upper bound of $o(\gamma^{w_\chi})$, we compute $\NN(\gamma^{w_\chi})=(\gamma^{w_\chi})^{w_\chi/v_\chi}$. This implies that $o(\gamma^{w_\chi})$ divides $w_\chi/v_\chi=:l^r$ and thus $s\leq w_\chi/v_\chi$.
To compute a lower bound, we first identify
$$\QQ^{\Q_l(\eta)}\Gamma^{w_\chi}\cong \Qu(\Z_l[\eta][[T]]),\ \ \gamma^{w_\chi}\leftrightarrow 1+T.$$
Now, assume that the order of $\gamma^{w_\chi}$ is less than $w_\chi/v_\chi$, i.e.~that there exists an $a\in \Qu(\Z_l[\eta][[T]])$ with $\NN(a)=(1+T)^{l^t}$ where $l^t<w_\chi/v_\chi=l^r$. 
 By definition, $a=\frac{f(T)}{g(T)}$ with $f(T),\ g(T)\in \Z_l[\eta][[T]]$. Let $\ell$ generate the maximal ideal of $\Z_l[\eta]$ above $l$. Then, the Weierstra\ss \ preparation theorem implies that
$$a=\frac{\ell^{n_1}F(T)u_1}{\ell^{n_2}G(T)u_2}$$
with $F(T)$ and $G(T)$ distinguished polynomials in $\Z_l[\eta][[T]]$ and $u_1$ and $u_2$ units in $(\Z_l[\eta][[T]])^\times$. The norm $\NN$ is now the Galois norm  of the field extension $\Q_l(\eta)/L$ and thus $\NN(\ell)$ generates the maximal ideal of the ring of integers $\o_L$ of $L$ because the extension is totally ramified. Furthermore, the norms of the distinguished polynomials resp.~units are distinguished polynomials resp.~units in $\o_L[[T]]$. Applying $\NN$ to $a$ thus yields
$$(1+T)^{l^t}=\NN(a)=\frac{\NN(\ell)^{n_1}\NN(F(T))\NN(u_1)}{\NN(\ell)^{n_2}\NN(G(T))\NN(u_2)}$$
 and therefore, with $u:=u_1/u_2\in (\Z_l[\eta][[T]])^\times$, we see
$$\NN(\ell)^{n_2}\NN(G(T))(1+T)^{l^t}=\NN(\ell)^{n_1}\NN(F(T))\NN(u).$$
Because $(1+T)^{l^t}$ is a unit, the Weierstra\ss \ preparation theorem for $\o_L[[T]]$ now implies that $n_1=n_2$, $\NN(F(T))=\NN(G(T))$, and thus that $F(T)$ and $G(T)$ only differ by a unit in $(\Z_l[\eta][[T]])^\times$. We conclude that $a\in (\Z_l[\eta][[T]])^\times$. 

We set $a=\sum_{i=0}^\infty a_i T^{i}$ with $a_i\in \Z_l[\eta]$ for all $i\geq 1$ and $a_0\in (\Z_l[\eta])^\times$, thus
\begin{align*}
(1+T)^{l^t}&=1+l^t T+...+T^{l^t}=\prod_{\sigma\in G(\Q_l(\eta)/L)}\left(\sum_{i=0}^\infty a_i T^{i}\right)^\sigma\\
&=\prod_{\sigma\in G(\Q_l(\eta)/L)}\left(\sum_{i=0}^\infty a_i^\sigma T^{i}\right)=\NN(a_0)+a_0 \Tr(a_1)T+...
\end{align*}
with $\Tr$ the trace of the field extension $\Q_l(\eta)/L$. Comparing the coefficients on both sides, we first see that $N(a_0)=1$ and therefore we can assume $a_0=1$ without loss of generality (otherwise divide $a$ by $a_0\neq 0$, this new $a$ has the same norm as the old one). Next, we consider the condition $l^t=\Tr(a_1)$. This equation is now to be read in $\Q_l(\eta)/L$. Set $\Q_l(\eta)=:\Q_l(\zeta)$ with $\zeta$ a  primitive $l$-power root of unity. We achieve $a_1=\alpha_0+\alpha_1 \zeta+...+\alpha_{l^r-1}\zeta^{l^r-1}$ with $\alpha_i\in \o_L$ for $0\leq i<l^r$, recall $l^r=[\Q_l(\eta):L]$. Thus 
\begin{align*}
l^t&=\Tr(a_1)=\Tr(\alpha_0+\alpha_1 \zeta+...+\alpha_{l^r-1}\zeta^{l^r-1})\\
&=l^r\alpha_0+\alpha_1\Tr(\zeta)+...+\alpha_{l^r-1}\Tr(\zeta^{l^r-1}).
\end{align*}
Next, we compute the trace of the powers of $\zeta$. For this, we look more closely at the field extension $\Q_l(\zeta)/L$. As $\Q_l(\zeta)$ is cyclic over $\Q_l$, we conclude that  $L$ itself is a cyclotomic field $L=\Q_l(\xi)$ with $\xi=\zeta^{l^r}$. Thus, the minimal polynomial of $\zeta$ over $L$ is $p(X)=X^{l^r}-\xi$. This implies that $\Tr_{\Q_l(\zeta)/L}(\zeta)=0$. Analogously, we obtain
$$ \Tr_{\Q_l(\zeta^{i})/L}(\zeta^i)=0   \mbox{ for all } 0\leq i<l^r$$
and furthermore, we set $i=l^\nu\iota$ with $0\leq \iota <l$ for all $0\leq i<l^r$. Obviously, we have $\Q_l(\zeta^i)=\Q_l(\zeta^{l^\nu})$. With this, we  compute
\begin{align*}
\Tr(\zeta^{i})&=\Tr_{\Q_l(\zeta)/L}(\zeta^{i})=\Tr_{\Q_l(\zeta^i)/L}\circ\Tr_{\Q_l(\zeta)/\Q_l(\zeta^i)}(\zeta^{i})\\
&=\Tr_{\Q_l(\zeta^{l^\nu})/L}(l^\nu\zeta^{i})=l^\nu\Tr_{\Q_l(\zeta^{i})/L}(\zeta^{i})=0.
\end{align*}
Therefore, we finally see
$$l^t=\Tr(a_1)=l^r\alpha_0,\  \mbox{i.e.} \ \alpha_0=l^{t-r}\in \Z_l[\eta].$$
This is not possible for $t<r$, and we thus have a contradiction to the assumption that the order of $\gamma^{w_\chi}$ was smaller than $w_\chi/v_\chi$.

We obtain $s\geq w_\chi/v_\chi$ and finally $s_D=s=w_\chi/v_\chi$. This shows both (iii) and (iv).

For the claim $(ii)$, we have seen that $\QQ^{\Q_l(\eta)}\Gamma^{w_\chi}$ is a maximal subfield of the cyclic skew field $D$ and thus is a splitting field of $D$ and $A$ (see e.g.~\cite[Thm (30.8)]{Re}).
\hfill$\Box$

\begin{corollary}
The $\QQ^c\Gamma_0$-algebra $\QQ^cG$ splits.
\end{corollary}

Finally, we give an example of a Galois group $G$ which causes nontrivial Schur indices\footnote{A first such example was given by A.~Weiss (unpublished).}:

Let $l=3$ and $H=\langle h\rangle$ be the cyclic group of order $9$. $G=H\rtimes \Gamma$ is determined by the action $h^{\gamma}=h^4$. 
Define the absolute irreducible constituent $\eta$ of $\res_G^H\chi$ by $\eta(h)=\zeta_9$ with $\zeta_9$ a primitive ninth root of unity, thus $St(\eta)=H\times\langle\gamma^3\rangle$ and $\chi=\ind_{St(\eta)}^G\chi'$ with $\chi'(h)=\eta(h)$ and $\chi'(\gamma^3)=1$. 

Then, $\chi(1)=3$, $w_\chi=3$ and $\Gal=G(\Q_l(\eta)/\Q_l)=G(\Q_l(\zeta_9)/\Q_l)=\langle\sigma\rangle $ with $\eta(h)^{\sigma}=\zeta_9^{\sigma}=\zeta_9^2$. Therefore, the action of $\gamma$ on $\eta(H)$ can be expressed as Galois action: $\eta(h)^\gamma=\eta(h)^4=\eta(h)^{\sigma^2}$ which implies $v_\chi=1$ and $G_0=\langle\sigma^2\rangle $.  

We achieve that the simple component $A$ of $\QQ G$ determined by $\chi$ has Schur index $s_D=w_\chi/v_\chi=3$ and  dimension $\dim_{Z(A)}A=\chi(1)^2=9=s_D^2$ and is therefore the cyclic skew field
$$D=(\QQ^{\Q_l(\zeta_9)}\Gamma^3/\QQ^{\Q_l(\zeta_3)}\Gamma^3, \sigma_{v_\chi}=\sigma^2, \gamma^3)$$
with $\zeta_3$ a primitive third root of unity. 

This example can be realized as Galois group over $\Q$ according to the situation of the Iwasawa extension examined earlier:
First, we show that $\Z/9\rtimes \Z/3$ can be realized as Galois group over $\Q$. With this realization, we then construct a field extension of $\Q$ with Galois group $G$.

 Clearly, $\Q(\zeta_9+\zeta_9^{-1})=\Q_{(1)}\subsetneq \Q_\infty$ is a subfield of the cyclotomic $\Z_3$-extension $\Q_\infty/\Q$ and $\Z/3\cong G(\Q(\zeta_9+\zeta_9^{-1})/\Q)$. For the whole group $\Z/9\rtimes\Z/3$, observe that in \cite{Ne} this is formulated as imbedding problem $\mathcal{E}_{\Z/3}(G(\Q^c/\Q), \Z/9\rtimes \Z/3)$
 corresponding to the diagram
 \[ \xymatrix{
        		 & {G(\Q^c/\Q)} \ar[d]^{\varphi}\\
        {\Z/9\rtimes \Z/3} \ar[r]^/2.5mm/{f} & {\Z/3} 
  						 } \]
 with  epimorphisms $f$ and $\varphi$. The problem is whether there exists a surjective homomorphism $\psi: G(\Q^c/\Q)\rightarrow \Z/9\rtimes \Z/3$ with $f\circ\psi=\varphi$. This is shown to be equivalent to the existence of a Galois extension $K_{(1)}\supseteq \Q(\zeta_9+\zeta_9^{-1})\supseteq \Q$ with Galois group $G(K_{(1)}/\Q)\cong\Z/9\rtimes \Z/3$ such that the canonical projection $G(K_{(1)}/\Q)\rightarrow G(\Q(\zeta_9+\zeta_9^{-1})/\Q)$ coincides with $\Z/9\rtimes \Z/3\rightarrow \Z/3$. 

\begin{floatingtable}[l]
	{
   \xymatrix@!C=0.7pc{
        		&	& & {K_\infty}\ar@{-}[d]\\
        		&&& {\Q_\infty}\ar@/^10mm/@{-}[dldldl]^{\cong\Gamma}\\
       		 & {\ K_{(1)}} \ar@{-}[-2,2] \ar@{-}[d]^{H}  \\
        {\ \ K_{(0)}} \ar@{-}[ur] \ar@{-}[d] &{\ \Q_{(1)}}\ar@{-}[-2,2]\\
        {\Q} \ar@{-}[ur]^{3}      
    }} 
\end{floatingtable}

 An affirmative answer to this question is given by \cite[Cor 6]{Ne}. To apply this, note that this corollary is formulated under the additional assumption that the kernel of $f$ is a pro-solvable group of finite exponent prime to the number $m(K_{(1)})$ of roots of unity in the field $K_{(1)}$ (see \cite[p. 157]{Ne}).  As $K_{(1)}$ is of odd index over $\Q(\zeta_9+\zeta_9^{-1})$, we know that $K_{(1)}$ is totally real and therefore $m(K_{(1)})=2$. Thus, $m(K_{(1)})$ is prime to the exponent  of the kernel of $f$. Now, we conclude with \cite[Cor 6]{Ne} that this imbedding problem  has a proper solution because we may choose a splitting group extension $\Z/9\rtimes\Z/3\rightarrow \Z/3$.

For the realization of $G=H\rtimes \Gamma\cong \Z/9\rtimes \Gamma$, consider the  diagram
with $\Q_{(1)}=\Q(\zeta_9+\zeta_9^{-1})$. By the above defined action  of $\Gamma$ on $H$, we see that $K_\infty$ is abelian over $\Q_{(1)}$ and therefore $G(K_\infty/\Q_\infty)\cong H$ and $G(K_\infty/\Q)\cong H\rtimes \Gamma$ as desired. Finally, note that $K_{(1)}$ is of odd index over $\Q$ and thus the constructed realization of $G$ is totally real.

\begin{remark}
\upshape{There is the natural question whether the structural results achieved in the present chapter can be generalized to non-pro-$l$ groups. However, already the split case $G=H\times \Gamma$ indicates that this is not possible. When $H$  is an $l$-group, the Schur indices of the Wedderburn components to all irreducible $\Q_l^c$-characters $\chi$ of $G$ with open kernel are trivial, because in this case $w_\chi=1$. But if $H$ is not an $l$-group, then the Schur index $s_D$ of the component corresponding  to  $\chi$ is essentially only restricted by $s_D\mid \chi(1)$ (compare e.g.~\cite{Lo}). }
\end{remark}

\section[Completion and cohomological dimension]{Completion and cohomological dimension}\label{completion}
\subsection{The completed algebra $\QQ_\wedge G$}

We again consider pro-$l$ Galois groups $G$ and use our results on the structure of $\QQ G$ to examine the completed algebra $\QQ_\wedge G$.
Let $A$ be the Wedderburn component of $\QQ G$ corresponding to the irreducible $\Q_l^c$-character $\chi$ and $D$ its underlying skew field.
The centre 
$$Z(D)=\QQ^{L}\Gamma^{w_\chi}=L\otimes_{\Q_l}\QQ \Gamma^{w_\chi}\cong L\otimes_{\Q_l}\Qu(\Z_l[[T]])$$
  is not complete with respect to any $\mathfrak{p}$-adic valuation where $\mathfrak{p}$ denotes the prime ideal $(\ell)$, $(T)$ or  $(f(T))$ for an irreducible distinguished polynomial of $\o_L[[T]]$ (with $\ell\in\o_L$ of value $1$). For example, consider for brevity $L=\Q_l$. Then, the sequence $a_n=\sum_{i=0}^n\frac{l^{i}}{T^{i}}$ is Cauchy with respect to the $l$-adic valuation but does not converge in $\Qu(\Z_l[[T]])\subsetneq \Q_l((T))$.

 For the examination of $\QQ_\wedge G$, we will need the concept of higher dimensional local fields.
 
\begin{definition}
A complete discrete valuation field $Q$ is called $n$-dimensional local field if there exists a chain of fields $Q=Q^{(n)},Q^{(n-1)},...,Q^{(0)}$ such that $Q^{(i+1)}$ is a complete discrete valuation field with residue field $Q^{(i)}$ and $Q^{(0)}$ is a finite field. 
\end{definition}

We first cite an example of higher dimensional local fields, the so-called standard fields (see \cite[p.~6]{Zh}).
For a complete discrete valuation field $F$ with residue field $\overline{F}$, we set
$$K_F:=F\{\{T\}\}:=\left\{\sum_{-\infty}^\infty a_iT^{i}:\ a_i\in F,\ \inf \{v_F(a_i)\}>-\infty,\ \lim_{i\rightarrow -\infty}v_F(a_i)=\infty\right\}.$$
It is a complete discrete valuation field with valuation $v_{K_F}(\sum a_iT^i)=\min\{v_F(a_i)\}$ and residue field $\overline{F}((\overline{T}))$. For a  local field $F$, the fields
$$F\{\{T_1\}\}...\{\{T_m\}\}((T_{m+2}))...((T_n)),\ \ 0\leq m\leq n-1,$$
are $n$-dimensional local fields. They are called standard fields.

Now, we can compute the completions of $\QQ^L\Gamma^{w_\chi}$.

\begin{lemma}\label{l4.1}
Let $X$ denote a transcendental element over $L$ and $\QQ_\wedge^L\Gamma^{w_\chi}$ the completion of  $\QQ^L\Gamma^{w_\chi}$ with respect to the $\frak{p}$-adic valuation. Then $\QQ_\wedge^L\Gamma^{w_\chi}$ is a two-dimensional local field and
\begin{enumerate}
\item $\QQ_\wedge^L\Gamma^{w_\chi}\cong L((X))$ if $\ \p=(T)$,
\item $\QQ_\wedge^L\Gamma^{w_\chi}\cong \left(L\otimes_{\Z_l}\Z_l[[T]]/(f)\right)((X))$ if $\ \p=(f)$ for an irreducible distinguished polynomial $f$,
\item $\QQ_\wedge^L\Gamma^{w_\chi}\cong L\{\{T\}\}$
 if $\ \p=(\ell)$.
\end{enumerate}
\end{lemma}
\textbf{Proof:}
The strategy of the proof is to compute the residue fields and then to choose the right case of the classification theorem in \cite{Zh} which classifies higher dimensional local fields up to isomorphism. 

We set $\Lambda^{\o_L}\Gamma^{w_\chi}\cong \o_L[[T]]=:R$ with $\xi$ a primitive $l$-power root of unity such that $\o_L=\Z_l[\xi]$ and $\QQ^L\Gamma^{w_\chi}\cong\Qu(R)=:Q$. Furthermore, let $(.)_{\b}$ resp.~$(.)_\wedge$ denote the localization with respect to the prime lying above $\p$ resp.~the $\p$-adic completion.

The first case $\p=(T)$ is a special case of the second case and we do not give an extra proof here.

In the second case $\p=(f)$ for a distinguished irreducible polynomial $f$, we see that the valuation ring is $\o:=(\o_L[[T]])_\b$ and therefore 
$$(Q_\wedge)^{(1)}=\o_\wedge/(f)=\o/(f)=(L\otimes_{\Z_l}\Z_l[[T]])/(f).$$
 Because $f$ is as irreducible distinguished polynomial either the monomial $T$ or Eisenstein, $(Q_\wedge)^{(1)}$ is totally ramified over $\Q_l$ and thus 
 $$(Q_\wedge)^{(0)}=\F_l.$$
  Hence, we conclude that $Q_\wedge$ is a $2$-dimensional local field and, moreover, that it is isomorphic to $\left((L\otimes_{\Z_l}\Z_l[[T]])/(f)\right)((X))$.

Finally, we show the third isomorphism for the prime $\p=(\ell)$. Because $l$ is totally ramified in $\Q_l(\xi)$, the prime ideal above $l$ is generated by, say,  $\ell=1-\xi$. With the valuation ring $\o=R_\bullet$, we compute 
$$(Q_\wedge)^{(1)}=\o_\wedge/(1-\xi)=\o/(1-\xi)=R_\b/(1-\xi)\supseteq R/(1-\xi)=\F_l[[\overline{T}]].$$
This is not yet a field but only the ring of integers of the local field $\F_l((\overline{T}))$. Hence, we achieve $(Q_\wedge)^{(1)}\supseteq \F_l((\overline{T}))$. By $(\Z_l[\xi])_\b=\o_L=\Z_l[\xi]$, the coefficient ring of $R$ does not change on passing to the localization $R_\bullet$  and therefore 
$$(Q_\wedge)^{(1)}=\F_l((\overline{T}))\ \ \mbox{and}\ \ (Q_\wedge)^{(0)}=\F_l[[\overline{T}]]/(\overline{T})=\F_l.$$ 
The classification theorem now implies that $Q_\wedge$ is a finite extension of a standard field $F\{\{X\}\}$ with local field $F$.

 For the computation of the exact type of $Q_\wedge$,  
 we follow the construction in the proof of the classification theorem:

First,  $\Q_l\{\{T\}\}$ is a  complete discrete valuation field of characteristic $0$ and $(l)$ is a prime ideal in $\o_{\Q_l\{\{T\}\}}$. The valuation is with respect to $l$, thus its absolute index of ramification is $e(\Q_l\{\{T\}\}):=v_{\Q_l\{\{T\}\}}(l)=1$. Furthermore, it has the same residue field $\overline{\Q_l\{\{T\}\}}=\F_l((\overline{T}))$ as $Q_\wedge$. Then \cite[II.5.6]{FV} implies that $Q_\wedge$ can be viewed as  a finite extension of $\Q_l\{\{T\}\}$.  

It remains to show that $Q_\wedge=L\{\{T\}\}$. To do so, we compare $LQ_\wedge=Q_\wedge$ and $L\Q_l\{\{T\}\}=L\{\{T\}\}$. For the ramification index, we have $$e(LQ_\wedge/L\Q_l\{\{T\}\})=e(Q_\wedge/L\{\{T\}\})=1$$
 because $1-\xi$ induces the discrete valuations of the two fields. As both fields are complete discrete valuation fields we conclude
$$[Q_\wedge:L\{\{T\}\}]=n=e(Q_\wedge/L\{\{T\}\})f(Q_\wedge/L\{\{T\}\})=[\overline{Q_\wedge}:\overline{L\{\{T\}\}}].$$
It therefore follows that the degree of the  field extension equals the residue class degree. Finally
\begin{align*}[Q_\wedge:L\{\{T\}\}] &=[\overline{Q_\wedge}:\overline{L\{\{T\}\}}]=[\F_l((\TT)):\F_l((\TT))]=1
\end{align*}
implies that 
$Q_\wedge=L\{\{T\}\}$
is a standard field.
\hfill$\Box$

 In particular, we get 
\begin{corollary}\label{corollary}
With the notations of Lemma \ref{l4.1}, we get for the residue fields
\begin{enumerate}
\item $\overline{\QQ_\wedge^L\Gamma^{w_\chi}}=\overline{\QQ^L\Gamma^{w_\chi}}\cong L$ if $\ \p=(T)$,
\item $\overline{\QQ_\wedge^L\Gamma^{w_\chi}}=\overline{\QQ^L\Gamma^{w_\chi}}\cong (L\otimes_{\Z_l}\Z_l[[T]])/(f)$ if $\ \p=(f)$ for an irreducible distinguished polynomial,
\item $\overline{\QQ_\wedge^L\Gamma^{w_\chi}}=\overline{\QQ^L\Gamma^{w_\chi}}\cong \F_l((\TT))$ if $\ \p=(\ell)$.
\end{enumerate}
\end{corollary}

\textbf{Proof:} We have shown this in the proof of Lemma \ref{l4.1}.\hfill$\Box$

\begin{corollary}
The cohomological dimension of the completed field is
$$\cd(\QQ^L_\wedge\Gamma^{w_\chi})=3.$$
\end{corollary}

\textbf{Proof:}
In the case $\p\neq (\ell)$, this follows directly from the fact that every local field has cohomological dimension $2$ and \cite[Thm (6.5.15)]{NSW}. 

For $\p=(\ell)$, the claim is stated in \cite[end of \S 3]{Mo}\footnote{A.~Weiss has drawn my attention to this paper.}.\hfill$\Box$

Let $D_\wedge=Z(D)_\wedge\otimes_{Z(D)}D$ be the completion of $D$ with respect to the $\p$-adic valuation. Then $D_\wedge$ is a central simple algebra over $Z(D_\wedge)=Z(D)_\wedge=\QQ_\wedge^L\Gamma^{w_\chi}$  and $[D_\wedge:Z(D_\wedge)]=[D:Z(D)]=s_D^2$ (see \cite[Hilfssatz 14.2]{Hu}).

\begin{proposition}
Let $D_\wedge$ be as above where the completion is with respect to $\p\neq (\ell)$. Then $$SK_1(D_\wedge)=1.$$
\end{proposition}

\textbf{Proof:}
Let $\overline{D_\wedge}$ denote the residue skew field of the completed skew field $D_\wedge$.
By the above, $\ch(\overline{\QQ_\wedge^L\Gamma^{w_\chi}})=0$ and therefore $\overline{\QQ_\wedge^L\Gamma^{w_\chi}}$ is perfect. In particular, the field extension $Z(\overline{D_\wedge})/\overline{\QQ_\wedge^L\Gamma^{w_\chi}}$ is separable. Moreover, we have seen that $\overline{\QQ_\wedge^L\Gamma^{w_\chi}}$ is a local field. Then the proposition follows readily from Korollar 7 in \cite{Dra}. We only have to substitute $k$ by $\QQ_\wedge^L\Gamma^{w_\chi}$ and Draxl's skew field $D$ is our $D_\wedge$. Observe that local fields are reasonable (vern\"unftig).\hfill$\Box$

 As these $D_\wedge$ are the centres of the Wedderburn components of $\QQ_\wedge G$, we even have proven 
\begin{corollary}
Let $\QQ_\wedge G$ be the Iwasawa algebra completed with respect to the prime $(T)$ or $(f(T))$ for an irreducible distinguished polynomial $f(T)$. Then
$$SK_1(\QQ_\wedge G)=1.$$
\end{corollary}

Now, we consider the case when $\ch(\overline{\QQ_\wedge^L\Gamma^{w_\chi}})=l$, i.e.~the completion  is with respect to the $(\ell)$-adic valuation. Then,
$\overline{\QQ_\wedge^L\Gamma^{w_\chi}} \cong \F_l((\overline{T}))$ is not perfect. Thus, $Z(\overline{D_\wedge})/\overline{\QQ_\wedge^L\Gamma^{w_\chi}}$ might not be separable and we are no longer in the situation of Draxl's Korollar 7. Yet, the following results can be transferred to this situation.

\begin{proposition}
Let $D_\wedge$ be as above where the completion is with respect to $\p = (\ell)$. Then, $D_\wedge$ is a skew field with Schur index $s_{D_\wedge}=s_D$ and its residue skew field $\overline{D_\wedge}$ is commutative.
\end{proposition}
\textbf{Proof:}
Let $M=\QQ^{\Q_l(\eta)}\Gamma^{w_\chi}=\QQ^{\Q_l(\zeta)}\Gamma^{w_\chi}$ be the maximal subfield of $$D=(\QQ^{\Q_l(\eta)}\Gamma^{w_\chi}/\QQ^{L}\Gamma^{w_\chi}, \sigma_{v_\chi}, \gamma^{w_\chi}),$$
 i.e.~$s_D=[D:M]=[M:Z(D)]$. The prime ideal $\p=(\ell)$ in $Z(D)=\QQ^{L}\Gamma^{w_\chi}$ is generated by $(1-\xi)$. In $M/Z(D)$, this prime ideal is totally ramified and undecomposed. Therefore, $M\otimes_{Z(D)}Z(D)_\wedge=M_\wedge$ is a field of degree $[M_\wedge:Z(D)_\wedge]=s_D$.
$M_\wedge$ is thus a maximal subfield in $D_\wedge$ because
$$Z(D)_\wedge\stackrel{s_D}{\subseteq}M\otimes_{Z(D)}Z(D)_\wedge\stackrel{s_D}{\subseteq} D\otimes_{Z(D)}Z(D)_\wedge=D_\wedge.$$
In particular, we conclude
$$D_\wedge=(M_\wedge/Z(D_\wedge), \sigma_{v_\chi}, \gamma^{w_\chi}).$$
Next, we show that this crossed product has Schur index $s_D$. We again have to check that the order $o(\gamma^{w_\chi})$ in $$Z(D_\wedge)^\times/\NN_{M_\wedge/(Z(D)_\wedge)}(M_\wedge)^\times=:Z(D_\wedge)^\times/\NN(M_\wedge)^\times$$
 is exactly $s_D$. Again, $o(\gamma^{w_\chi})$ divides $s_D=l^r$ because of $\NN(\gamma^{w_\chi})=(\gamma^{w_\chi})^{s_D}$. 
 
 Now, we assume that $o(\gamma^{w_\chi})=l^t$ with $t<r$. Then, there exists an $a\in M_\wedge$ such that $\NN(a)=(\gamma^{w_\chi})^{l^t}$. Therefore, $a$ is integral as $\gamma^{w_\chi}$ is. Furthermore, the residue fields of $M_\wedge$ and $Z(D)_\wedge$ coincide as $\p$ is totally ramified, i.e.~$G(\overline{M_\wedge}/\overline{Z(D_\wedge)})=\langle\overline{\sigma_{v_\chi}}\rangle=1$. We achieve
$$(\overline{\gamma^{w_\chi}})^{l^t}=\overline{\NN(a)}=\prod_{j=0}^{s_D-1}\overline{a}^{\overline{\sigma_{v_\chi}}^j}=\overline{a}^{s_D}=\overline{a}^{l^r}.$$
If, as we assume, $t<r$, then $\overline{\gamma^{w_\chi}}$ is an $l$-th power. To show that this is not possible, we use the isomorphisms $\QQ^{L}_\wedge\Gamma^{w_\chi}\cong L\{\{T\}\}$ of Lemma \ref{l4.1} and $\overline{\QQ^{L}_\wedge\Gamma^{w_\chi}}\cong \F_l((\TT))$  of Corollary \ref{corollary}.
Indeed, 
\begin{align*}
\overline{\gamma^{w_\chi}}\leftrightarrow 1+\TT=\left(\sum_{i=-n}^\infty \alpha_i \TT^i\right)^l=\sum_{i=-n}^\infty \alpha_i^l\TT^{il}
\end{align*}
is a contradiction.

Finally, we show that $\overline{D_\wedge}$ is commutative. We use valuation theory for the skew field $D_\wedge$, for details see e.g.~\cite[\S 3]{Dra}.

 We start with computing $[\overline{D_\wedge}:\overline{Z(D_\wedge})]$. For this, let $v$ denote the $\p$-adic valuation of $Z(D_\wedge)$ induced by $l$, i.e.~$\p=(1-\xi)$. The extension $w$ of $v$ to $D_\wedge$ is defined by
$$w(d)=\frac{1}{s_{D_\wedge}}v(\nr_{D_\wedge/Z(D_\wedge)}(d))$$
for every $d\in D_\wedge$.  We then compute the ramification index $$e(D_\wedge/Z(D_\wedge))=[w(D_\wedge^\times):v(Z(D_\wedge)^\times )].$$
 The definition of $w$ implies that $e(D/Z(D))\leq s_{D_\wedge}$. As $\Q_l(\eta)=\Q_l(\zeta)$ is totally ramified over $L=\Q_l(\xi)$, we have $$\NN_{\Q_l(\zeta)/\Q_l(\xi)}(1-\zeta)=\nu\cdot(1-\xi)$$
 for a unit $\nu$ in the valuation ring of $\Q_l(\xi)$.

 We now choose $d=1-\zeta\in M_\wedge^\times \subseteq D_\wedge^\times$ and compute
\begin{align*}
w(1-\zeta)&=\frac{1}{s_{D_\wedge}}v(\nr_{D_\wedge/Z(D_\wedge)}(1-\zeta))=\frac{1}{s_{D_\wedge}}v(\NN_{M_\wedge/Z(D_\wedge)}(1-\zeta))\\ &=\frac{1}{s_{D_\wedge}}v(\NN_{\Q_l(\zeta)/\Q_l(\xi)}(1-\zeta))
=\frac{1}{s_{D_\wedge}}v(\nu(1-\xi))=\frac{1}{s_{D_\wedge}}.
\end{align*}

This implies that $e(D_\wedge/Z(D_\wedge))\geq s_{D_\wedge}$ and finally $e(D_\wedge/Z(D_\wedge))=s_{D_\wedge}$. By the equation 
$$s_{D_\wedge}^2=[D_\wedge:Z(D_\wedge)]=e(D_\wedge/Z(D_\wedge))[\overline{D_\wedge}:\overline{Z(D_\wedge)}]=s_{D_\wedge}[\overline{D_\wedge}:\overline{Z(D_\wedge)}],$$
we achieve that
$$[\overline{D_\wedge}:\overline{Z(D_\wedge)}]=s_{D_\wedge}.$$

Next, we consider $N_\wedge:=Z(D_\wedge)(\gamma^{v_\chi})$. We have seen that   $Z(D_\wedge)\cong L\{\{T\}\}$ with $\gamma^{w_\chi}\leftrightarrow 1+T$. As
$\gamma^{v_\chi}$ commutes with $L$ and $\gamma^{w_\chi}$, this implies that $N_\wedge$ is commutative. Moreover,  $N_\wedge\cong L\{\{T\}\}[X]/(f(X))$ with  $f(X)=X^{w_\chi/v_\chi}-(1+T)=X^{s_{D_\wedge}}-(1+T)$. The polynomial $f(X)$ is irreducible modulo $\p$ and thus $f(X)$ is irreducible itself. 
Thus, $N_\wedge$ is a subfield $Z(D_\wedge)\subseteq N_\wedge\subseteq D_\wedge$ with $$[\overline{N_\wedge}:\overline{Z(D_\wedge)}]=[N_\wedge:Z(D_\wedge)]=s_{D_\wedge}=[\overline{D_\wedge}:\overline{Z(D_\wedge)}].$$
 We finally conclude that $\overline{D_\wedge}=\overline{N_\wedge}$ is commutative.
\hfill$\Box$

\subsection{Adapting Saito's result}
We show, using a geometrical result of S.~Saito, that the centres of the Wedderburn components of $\QQ G$ are of cohomological dimension 3. This reduces the uniqueness of $\Theta$ in the main conjecture of equivariant Iwasawa theory to a conjecture by Suslin which says that $SK_1(D)=1$ if $Z(D)$ is of cohomological dimension $3$ and $D$ is of finite degree over $Z(D)$. For details on the conjecture, we refer to \cite{Sus}. 
Our computations also justify the remark in \cite[p.~565]{Iw2} that the centre fields in $Z(\QQ G)$ have cohomological dimension $3$ by a result of Kato (which was published by S.~Saito in \cite{Sa}). 
 
\begin{theorem}\label{t4.1}
Let $D$ be the underlying skew field of a simple component of $\QQ G$. Then
$$\cd(Z(D))=3.$$
\end{theorem}

For the proof of the theorem, we essentially need Theorem 5.1 of \cite{Sa}:
\begin{theorem}[S.~Saito]\label{Saito}
Let $A$ be a $2$-dimensional excellent normal henselian local ring, $K_A$ its field of fractions and $F_A$ its algebraically closed residue field. Then
$$\cd_p(K_A)=2 \mbox{ for every prime number } p\neq \ch(K_A).$$
\end{theorem}
\textbf{Proof of Theorem \ref{t4.1}:}
We have $Z(D)=\QQ^L \Gamma^{w_\chi}=\QQ^{\Q_l(\xi)} \Gamma^{w_\chi}$ for a primitive $l$-power root of unity $\xi$. In the completed situation, we have seen that the cohomological dimension is $3$. Thus, in the uncompleted case, the cohomological dimension is at least $3$ by 
the following:
For any field $E$, let $\Gal(E)$ denote its absolute Galois group. Then, $\Gal(\QQ^{L}_\wedge \Gamma^{w_\chi})$ is a closed subgroup of $\Gal(\QQ^{L} \Gamma^{w_\chi})$ because it is the decomposition group of some prime above $\p$, where the completion is with respect to $\p$. Thus, 
\cite[(3.3.5)]{NSW} shows the claim:
$$3= \cd(\QQ_\wedge^L\Gamma^{w_\chi})\leq\cd(\QQ^L\Gamma^{w_\chi}).$$

Furthermore, we see that $\cd(\QQ^{\Q_l(\xi)} \Gamma^{w_\chi})\leq \cd(\QQ \Gamma^{w_\chi})$ again by \cite[~(3.3.5)]{NSW} because $\Gal(\QQ^{\Q_l(\xi)} \Gamma^{w_\chi})$ is a closed subgroup of $\Gal(\QQ \Gamma^{w_\chi})$. 

We thus consider the problem whether $\cd(K_A)=3$ for $K_A=\QQ\Gamma^{w_\chi}$ and $A=\Lambda \Gamma^{w_\chi}$. Then, $\cd(\QQ^L \Gamma^{w\chi})\leq 3$ follows.

The residue field $F_A=A/\m\cong \Z_l[[T]]/(l,T)=\F_l$ is not algebraically closed and we thus can not apply Saito's theorem directly. To achieve the algebraic closure, we replace $A$ by the ring $A'=\varinjlim \Z_l[\zeta_i][[T]]=\bigcup\Z_l[\zeta_i][[T]]$ with $(l,i)=1$ and $\zeta_i$ a primitive root of order $i$.

 First, we look more closely to $A'$. Set $I=\N\verb=\=l\N$ the set of natural numbers coprime to $l$. This is filtered with the relation $i\leq j\ :\Leftrightarrow\ i\mid j$. We set $$A_i=\Z_l[\zeta_i][[T]]=\Z_l[[T]][\zeta_i]$$
 
 and define for $i\leq j$ the inclusion $\varphi_{ij}:A_i\rightarrow A_j$. Then, $A':=\varinjlim A_i=\bigcup A_i$ is the filtered inductive limit of $(A_i,\varphi_{ij})$.
 
Moreover, $A'$ is locale-ind-\'{e}tale in the sense of \cite[p.~80]{Ray} by the following: locale-ind-\'{e}tale means that $A'$ is the filtered inductive limit of some $(A_i,\varphi_{ij})$, where the $\varphi_{ij}$ are local morphisms and the  $A_i$ are locale-\'{e}tale $A$-algebras. 
First, $A_i$ is local with maximal ideal $\m_i=(l,T)$ because $l$ and $i$ are coprime. 
Next, $\varphi_{ij}: A_i\rightarrow A_j$ is local, i.e. $\varphi_{ij}^{-1}(\m_j)\subseteq \m_i$, because $\m_i$ and $\m_j$ are both generated by $l$ and $T$.
 $A_i$ is said to be locale-\'{e}tale over $A$ if it is the localization $B_\n$ of an \'{e}tale $A$-algebra $B$ by a prime ideal $\n$ over $\m$. We choose $B=A_i$ and $\n=\m_i=(l,T)$. Then $B_\n=B=A_i$.
 Thus, we only have to show that $A_i$ is \'{e}tale over $A$. To do so, we use \cite[(18.4.5)]{EGA}. $A$ and $A_i$ are local rings. Moreover, $A_i=A[\zeta_i]=A[X]/(F(X))$, with $F(X)$ an irreducible polynomial dividing the $i$-th cyclotomic polynomial, is a finite $A$-algebra of finite presentation. $F$ is unitary and  separable in the meaning of \cite[p.~118]{EGA}, i.e. $F'(\zeta_i)\notin \m_i=(l,T)$ for the generator $\zeta_i$  of $A_i$ over $A$ with minimal polynomial $F$, which can be seen as follows: 
$F '(\zeta_i)\notin \m_i=(l,T)$ is equivalent to $F'(\zeta_i)\not\equiv 0\mod \m_i$ or $\overline{F'(\zeta_i)}\neq 0$ in $A_i/\m_i=\F_l(\overline{\zeta_i})=\Z_l[\zeta_i]/(l)$. Because $\Z_l[\zeta_i]$ is unramified over $\Z_l$, we know that $\overline{F}$ is also the minimal polynomial of $\overline{\zeta_i}$ over $\F_l$ and in particular irreducible and separable. Thus, $\overline{F'(\zeta_i)}\neq 0$ in $\F_l(\overline{\zeta_i})$ and hence $F$ is a unitary separable polynomial as needed.
Hence, \cite[(18.4.5)]{EGA} implies that $A_i$ is an \'{e}tale $A$-algebra.
 
Now, we have $K_{A'}=\Qu(A')=\Q_l^{\ur}\otimes_{\Q_l}K_A$ and  $$G(K_{A'}/K_A)\cong G(\Q_l^{\ur}/\Q_l)=\hat{\Z}.$$

Assume for the moment that $A'$ fulfils  the conditions of Saito's theorem 
and thus $\cd_p(K_{A'})=2$ for all $p\neq \ch(K_{A'})=0$, i.e.~$\cd(K_{A'})=2$. 
 
$\Gal(K_{A'})$ is a closed normal subgroup of $\Gal(K_{A})$, therefore $\cd(\hat{\Z})=1$ (see \cite[p.~140]{NSW}) and \cite[(3.3.7)]{NSW}   imply
\begin{align*}
\cd(K_{A})&\leq \cd(K_{A'})+\cd(\Gal(K_A)/\Gal(K_{A'}))\\
&=\cd(K_{A'})+\cd(G(K_{A'}/K_{A}))=\cd(K_{A'})+\cd(\hat{\Z})=3.
\end{align*}

Now, we show that $A'$ fits into the situation of Saito. First, we observe that by construction of $A'$, the residue field $F_{A'}$ of $A'$ is the separable closure of $F_A=\F_l$. 
Moreover, $F_A$ is perfect as a finite field  and thus its separable closure $F_{A'}$ is already algebraically closed, i.e.~$F_{A'}=\F_l^c$.

 Next, we show that $A'$ is local.
 We use \cite[Prop 1, p. 6]{Ray} which says that the filtered inductive limit $A'$ is local because the $\varphi_{ij}$ are local morphisms and the $A_i$ local rings. Moreover, the proof of \cite[Prop 1, p.~6]{Ray} shows that the maximal ideal of $A'$ is $\m'=(l,T)$.
 
Analogously, we show that $A'$ is henselian. The $A_i$ are henselian because they are local and complete in the $\m_i$-adic topology (see \cite[p.~242]{NSW}). Thus, $A'$ is henselian again by \cite[Prop 1, p.~6]{Ray}.
 
For $A'$ to be normal, we show that every integral quotient $\frac{a}{b}\in K_{A'}$ lies in $A'$. There exists a finite subextension $A\subseteq A_i\subsetneq A'$ with $a, b\in A_i$. Then, $A_i$ is normal as finite extension of the normal ring $A$. Hence, we conclude $\frac{a}{b}\in A_i$ and thus $\frac{a}{b}\in A'$. 

The Krull dimension of $A'$ is $2$ as needed: Let $(l,T)\supsetneq \p_1\supsetneq \p_2$ be a descending chain of prime ideals in $A'$. We cut down this chain to every $\Z_l[\zeta_i][[T]]=A_i$. Because the maximal ideal of this $\Z_l[\zeta_i][[T]]$ is also generated by $l$ and $T$, the section of $(l,T)$ remains the maximal ideal and $A_i\cap (l,T)\supsetneq A_i\cap \p_1$ (otherwise $\p_1=(l,T)$, a contradiction). We show that $\p_2=0$. For this, assume that $\p_2\neq 0$, i.e.~there exists an element $0\neq a\in \p_2$ and a certain $\zeta_i$ such that $a\in \Z_l[\zeta_i][[T]]$. But $\Z_l[\zeta_i][[T]]$ has dimension $2$ and therefore  $\p_1$ and $\p_2$ coincide for $\Z_l[\zeta_i][[T]]$. Thus, $\p_1$ and $\p_2$ coincide for all  $\Z_l[\zeta_j][[T]]$ with $\zeta_j$ a primitive root of unity of order prime to $l$ such that $\zeta_j\mid \zeta_i$. Finally, $\p_1=\p_2$ in $A'$ because these $\Z_l[\zeta_j][[T]]$ already generate $A'$.

It remains to show that $A'$ is excellent. 
 First, we check that $A'$ is Noetherian.
For this, we use \cite[Thm 3.3, p. 94]{Ray} which says that the locale-ind-\'{e}tale $A$-algebra $A'$ is Noetherian if and only if $A$ is Noetherian. But $A=\Z_l[[T]]$ is Noetherian because it is a power series ring over the Noetherian ring $\Z_l$ by \cite[Thm 9.4, p. 210]{La}.
Next, $A$ is a Noetherian local ring and complete in its $\m$-adic topology (see \cite[p.~242]{NSW}). Thus, $A$ is excellent by \cite[p. 260]{Ma}.
Then, \cite[Thm 5.3.iv]{Gr}\footnote{M. Nieper-Wi\ss kirchen has drawn my attention to this paper.}  implies that $A'$ is also excellent because the inclusion $A\rightarrow A'$ is ind-\'{e}tale and thus AF, compare \cite[Def 5.1.]{Gr}.

This concludes the proof of the theorem.
\hfill$\Box$

\begin{corollary}
For the group $\Gamma\cong \Z_l$ we have
$$\cd(\QQ \Gamma)=3.$$
\end{corollary}

\bibliographystyle{spmpsci} 
\bibliography{Literature}

\end{document}